\documentclass[12pt]{article}
\usepackage{amsmath,amssymb,amsthm,mathtools,enumitem}
\usepackage{authblk}
\usepackage{color}
\usepackage{geometry}
\newtheorem{Theorem}{Theorem}[section]
\newtheorem{Lemma}{Lemma}[section]

\newtheorem{Definition}{Definition}[section]

\begin{document}

\def\eR{\mathbf{R}}
\def\Rd{{\eR}^d}
\def\Rdd{{\eR}^{d\times d}}
\def\Rdsym{{\eR}^{d\times d}_{sym}}
\def\eN{\mathbf{N}}
\def\eZ{\mathbf{Z}}
\def\LND{L^2_{\dvr}}
\def\CrlE{\mathcal{E}}
\def\CrlH{\mathcal{H}}
\def\CrlG{\mathcal{G}}
\def\IdM{\operatorname{I}}
\def\IdOp{\operatorname{Id}}
\def\LTD{\mathcal L^2_{n,\dvr}}
\def\LTDD{\dot{\mathcal L}^2_{\dvr}}
\def\WN{W^{1,2}_{n}}
\def\WND{W^{1,2}_{0,\dvr}}
\def\WDD{\dot{W}^{1,2}}
\def\WNDT{W^{1,2}_{n,\dvr}}
\def\dist{\operatorname{dist}}
\def\dd{\mbox{d}}
\def\VS{\mathcal{V}}
\def\WS{\mathcal{W}}
\def\Meas{\mathcal{M}}
\def\NnMeas{\mathcal{M}^+}
\newcommand{\TR}{\operatorname{Tr}}
\newcommand{\essinf}{\operatorname{ess\,inf}}
\newcommand{\esssup}{\operatorname{ess\,sup}}
\newcommand{\supp}{\operatorname{supp}}
\newcommand{\spn}{\operatorname{span}}
\newcommand\dx{\; \mathrm{d}x}
\newcommand\dy{\; \mathrm{d}y}
\newcommand\dz{\; \mathrm{d}z}
\newcommand\dt{\; \mathrm{d}t}
\newcommand\ds{\; \mathrm{d}s}
\newcommand\diff{\mathrm{d}}
\newcommand\dvr{\mathop{\mathrm{div}}\nolimits}
\newcommand\diam{\mathrm{diam}}
\newcommand\tder{\partial_t}
\newcommand\tr{\mathop{\mathrm{tr}}\nolimits}
\newcommand\curlP{\mathop{\mathcal{P}}\nolimits}
\newcommand\curlQ{\mathop{\mathcal{Q}}\nolimits}
\newcommand\lspan{\mathop{\mathrm{span}}\nolimits}
\newcommand{\IED}{IED}
\newcommand\Hext{H_{ext}}
\newcommand\colorstart{\color{black}} % change to some color for the revised version

\title{Dissipative solutions to a system for the flow of magnetoviscoelastic materials}
\author{Martin Kalousek\footnote{Institute of Mathematics, University of W{\"u}rzburg, Emil-Fischer-Str.\ 40, 97074 W{\"u}rzburg, Germany, Email: martin.kalousek@mathematik.uni-wuerzburg.de}\; and Anja Schl{\"o}merkemper\footnote{Institute of Mathematics, University of W{\"u}rzburg, Emil-Fischer-Str.\ 40, 97074 W{\"u}rzburg, Germany, Email: anja.schloemerkemper@mathematik.uni-wuerzburg.de}}
\date{\today}
\affil{}

\maketitle
\sloppy


\begin{abstract}
 We address the question of global in time existence of solutions to a magnetoviscoelastic system with general initial data. We show that the notion of dissipative solutions allows to prove such an existence in two and three dimensions. This extends an earlier result for the viscoelastic subsystem to the setting which includes the magnetization vector and its evolution in terms of a Landau-Lifshitz-Gilbert equation. 
  \end{abstract}

\section{Introduction}
Let $T>0$ and $\Omega\subset\eR^d$ a bounded domain with $d=2,3$. We analyze the following system of partial differential equations 
\begin{equation}\label{ClassForm}
\begin{aligned}
		\tder u+ (u \cdot \nabla) u -\Delta u+\nabla p&=\dvr(\nabla_FW(F)F^\top-\nabla M\odot \nabla M)+\nabla\Hext^\top M\\
		\dvr u&=0,\\
		\tder F+\dvr(u\otimes F)&=\nabla u F,\\
		\tder M+ (u\cdot \nabla) M +M\times(\Delta M+\Hext)&=\Delta M+M|\nabla M|^2-M(M\cdot\Hext)+\Hext,
	\end{aligned}
\end{equation}
where $u\!:(0,T)\times\Omega\to\eR^d$ denotes the velocity of the fluid, $p\!:(0,T)\times\Omega\to \eR$ the pressure, $F\!:(0,T)\times\Omega\to\eR^{d\times d}$ the deformation gradient \colorstart satisfying the constraint $\dvr F=0$\ in $(0,T)\times\Omega$, \color{black} and $M\!:(0,T)\times\Omega\to\eR^3$ the magnetization which additionally fulfills $|M|\equiv 1$ \ in $(0,T)\times\Omega$. For the assumptions on the elastic energy $W$ and the external magnetic field $\Hext$ we refer to \eqref{W1} and \eqref{HExtReg}, resp. We assume homogeneous Dirichlet boundary conditions for the velocity and homogeneous Neumann boundary condition for the magnetization, i.e.,
\begin{equation*}
	u=0,\ (\nabla M) n=0 \text{ on }(0,T)\times\partial\Omega.
\end{equation*}
The initial data
\begin{equation*}
	u(0,\cdot)=u_0,\ F(0,\cdot)=F_0,\ M(0,\cdot)=M_0 \text{ in }\Omega
\end{equation*}
with $\dvr u_0=0$, $\dvr F_0=0$, $|M_0|\equiv 1$ in $\Omega$ supplement system~\eqref{ClassForm}.

This system was derived via the energetic variational principle, cf.\ \cite{BFLS18,Fo16}. It models the evolution of magnetoviscoelastic materials that belong to a wider class of smart materials. They are characterized by the ability to change significantly, but in a controllable fashion, their mechanical properties under an external magnetic field. \colorstart Here, we have an additional constraint, namely $\dvr F=0$, which is needed for analytical reasons, cf.~\eqref{analyticalreason}.
However, the constraint is natural if one, for instance, has a constant matrix like the identity as an initial condition for $F$ since then Lemma~\ref{Lem:TranspEqRegData} yields this constraint automatically for sufficiently regular $F$ and $u$, see also \cite{LLZ05}. \color{black}

The goal of the paper consists in proving the global in time existence of a solution to system \eqref{ClassForm} with general initial data in two and three dimensions. This task is highly nontrivial also in the viscoelastic case, i.e., if $M\equiv 0$ and the system couples the momentum equation and the transport equation for the deformation gradient. Up to the authors' knowledge, the global in time existence of a weak solution for general initial data is still an open problem. The difficulty lies in the fact that the energy bound on $F$ yields only the compactness with respect to a weak topology, which does not allow to pass to the limit in the nonlinearity $\nabla_FW(F)F^\top$ in the stress term of the momentum equation \eqref{ClassForm}$_1$.
The coupling of the momentum equation and the Landau-Lifshitz-Gilbert equation in \eqref{ClassForm}$_4$ represents another difficulty for the analysis. Namely, the energy bound on $M$ provides only an $L^1$ bound with respect to the space variable of the term $\nabla M\odot\nabla M$ in the momentum equation. Hence, when one considers sequences approximating a solution, it is not possible to pass to the limit in the above mentioned nonlinear terms, at least in the sense of distributions.

Available existence results for systems with similar couplings either for a boundary value or a Cauchy problem were obtained under the assumption of a suitable closeness of initial data to the equilibrium possibly in a combination with the regularization of the transport equation for $F$, see e.g.\ \cite{LLZ08, LZ08,LLZ05} for the viscoelastic system and \cite{BFLS18,JLL19,KKS19} for the magnetoviscoelastic system. The articles \cite{JLL19,KKS19,LLZ05} also treat the local in time existence of a solution. 

In this article we address the question: Is there a notion of solution to \eqref{ClassForm} for which we can show the existence globally in time for general initial data? We consider the notion of dissipative solution that was introduced in the context of the incompressible Euler system in \cite{Li96} and later adopted for a hyperbolic system in \cite{FeRoScZa18}. Taking as an inspiration the existence result for a dissipative solution to the corresponding viscoelastic system \cite{K19}, we introduce a dissipative solution also for system \eqref{ClassForm} and prove its global in time existence (Theorem~\ref{Thm:Exis}). Roughly explained, the dissipative solution satisfies \eqref{ClassForm} in the sense of integral identities with
suitably regular test functions whose part corresponding to the right hand side of \eqref{ClassForm}$_1$ contains an extra
term regarded as a defect measure, cf.\ \eqref{WeakForms}$_1$. Moreover, a function called dissipation defect appears in the
energy inequality. This dissipation defect is attributed to singularities that may hypothetically
emerge during the fluid evolution. It dominates in a certain sense the additional term on the right
hand side of the integral formulation of \eqref{ClassForm}$_1$, see \eqref{DefMesEst}. Since the notion of a dissipative solution is quite weak, it is natural to study relations to other notions of a solution to \eqref{ClassForm} for which at least the local in time existence can be shown, cf., e.g., the dissipative-strong uniqueness proven in \cite{K19}. We plan to discuss these issues in a separate paper. % in order to keep a reasonable length of this one.

The strategy of the existence proof is to first consider an approximative system, see  \eqref{ApproxSysClass}, which has the same initial data for $u$ and $M$ and a regularized initial condition for $F$. Further it approximates the magnetization vector $M$ that needs to fulfil $|M|\equiv 1$ by $M^\varepsilon$ solving a parabolic equation with a penalizing term $\varepsilon^{-1}M^\varepsilon(|M^\varepsilon|^2-1)$. For showing existence of a weak solution to the approximative system (Lemma~\ref{Lem:ApproxEx}), we adopt ideas from \cite{AlSo,BMPS18}. In \cite{BMPS18}, the proof of existence of a solution to a problem with a coupling similar to the coupling of the momentum equation and the equation governing the evolution of $M$ is given. This is combined with techniques used in \cite{AlSo} to show existence of a weak solution to the Landau-Lifshitz-Gilbert equation without the convective term and an external magnetic field. The existence proof is based on a multi-level Galerkin scheme, which involves further parabolic regularizations. In particular, the challenging task here is the derivation of the energy inequality that turns out to be possible although approximations of $M^\varepsilon$ lack square integrable second derivatives up to the boundary that is assumed to be Lipschitz in our situation, see the end of Section~\ref{Sec:Stat} for a more detailed outline. The obtained energy inequality yields uniform bounds on a solution to the approximative problem that allow to let the regularizing parameter tend to zero and thus to obtain the existence of a dissipative solution.

The outline of the paper is as follows. In Section~\ref{Sec:Stat}, we first fix some notation and give a definition of dissipative solutions to system~\eqref{ClassForm}. In Theorem~\ref{Thm:Exis} we state the main result of the paper, which is the global in time existence of a dissipative solution for general initial data. The proof with all its steps is given in Section~\ref{Sec:DisExist}.
The appendix contains two technical lemmas. The first one is devoted to the equivalent formulations of the Landau-Lifshitz-Gilbert equation provided that its solution possesses sufficient regularity. The second lemma summarizes several assertions concerning the transport equation with regular data.

\section{Formulation of the results}\label{Sec:Stat}
We start by fixing some notation. The ball centered at $x$ with radius $r$ is denoted by $B(x,r)$. The centered dot $\cdot$ denotes the scalar product between vectors and matrices, respectively. If $a\in\eR^l$ and $B\in\eR^{m\times n}$, the outer product $a\otimes B$ denotes the tensor with components $a_iB_{jk}$, $i=1,\ldots, l$, $j=1,\ldots, m$, $k=1,\ldots,n$. \colorstart Further, $\nabla M \odot \nabla M$ is shorthand for $(\nabla M)^\top \nabla M$. \color{black} Generic constants are denoted by $c$. Let $\Omega\subset\eR^d$ be open. For $t>0$ and $\Omega\subset\eR^d$ we use the notation $Q_t$ for the time-space cylinder $(0,t)\times \Omega$. 
The space of Radon measures on $\overline\Omega$ is denoted by $\mathcal{M}(\overline\Omega)$ \colorstart with norm $\|\cdot \|_{\mathcal{M}}$, \color{black} the space of nonnegative Radon measures by $\mathcal{M}^+(\overline\Omega)$. 

Let $k\in\eN$ and $q\in[1,\infty]$. Then the standard Lebesgue and Sobolev spaces are denoted by $(L^q(\Omega),\|\cdot\|_q)$ and $(W^{k,q}(\Omega),\|\cdot\|_{k,q})$. For any Banach space  $X$ of scalar functions, we write $X^m$ for the corresponding space of  vector-valued functions with $m$ components each belonging to $X$. Similarly,  $X^{m\times n}$ denotes a Banach space of matrix-valued functions. In the notation of the corresponding norms, we often depress the dimension of the target space; we write e.g.\ $\|\cdot \|_{L^q(\Omega)}$ instead of $\|\cdot \|_{L^q(\Omega)^m}$, $ \|\cdot \|_{W^{k,q}(\Omega)}$ instead of $\|\cdot \|_{W^{k,q}(\Omega)^{m\times n}}$, and $\|\cdot \|_{L^q(0,t;L^r(\Omega))}$ instead of $\|\cdot \|_{L^q(0,t;L^r(\Omega)^{m})}$, etc. For the sake of clarity, the notation $(u,v)=\int_\Omega u(x)\cdot v(x)\dx$ is used for the scalar product in $L^2$.

If $X,Y$ are Banach spaces, the notation $X\hookrightarrow Y$ and  $X\stackrel{C}{\hookrightarrow} Y$ is used for expressing an embedding of $X$ to $Y$ that is continuous and compact, respectively. We denote the dual space of $X$ by $X^*$ and the corresponding duality pairing by $\left\langle\cdot,\cdot\right\rangle$. Furthermore, we need the set of weakly continuous mappings $C_w([0,T];X)$ that contains functions $f\in L^\infty(0,T;X)$ for which the real-valued mapping $t\mapsto\left\langle \phi, f(t)\right\rangle$ is continuous on $[0,T]$ for any $\phi\in X^*$. Further, we set 
\begin{equation*}
\begin{split}
		\LND(\Omega)&=\overline{\{v\in C^\infty_c(\Omega)^d:\ \dvr v=0\text{ in }\Omega\}}^{\|\cdot\|_{L^2}},\\
		\CrlH(\Omega)&=\{\Phi\in L^2(\Omega)^{d\times d}:\ \dvr \Phi=0\text{ in }\Omega\},
		\end{split}
\end{equation*}
where the distributional divergence of a $d\times d$ matrix-valued function $\Phi$ satisfies 
\begin{equation*}
	\left\langle \dvr\Phi,\varphi\right\rangle=-\int_\Omega \Phi\cdot(\nabla\varphi)^\top,\text{ for all }\varphi\in C^\infty_c(\Omega)^d.
\end{equation*} 
If $\Phi$ is smooth, we set $(\dvr(\Phi))_j=\sum_{i=1}^d\partial_i\Phi_{ij}$. That is, we here take the divergence of columns; if one defines the divergence of a tensor by taking the divergence of its rows, \eqref{ClassForm}$_4$ would read $\dvr F^\top = 0$. However, the same results of this article would hold true.

The following subspaces of $W^{1,2}(\Omega)^d$, $W^{1,2}(\Omega)^{d\times d}$ respectively are of relevance in this article: 
\begin{align*}
	\WS(\Omega)&=\{\Phi\in W^{1,2}(\Omega)^{d\times d}:\dvr \Phi=0\text{ in }\Omega\},\\
	\WND(\Omega)^d&=\overline{\{v\in C^\infty_c(\Omega)^d: \dvr v =0\text{ in }\Omega\}}^{\|\cdot\|_{W^{1,2}}}.
\end{align*}
Further, we use also the space $\VS(\Omega)$ defined as
\begin{equation*}
	\VS(\Omega)=\overline{\{v\in C^\infty_c(\Omega)^d: \dvr v =0\text{ in }\Omega\}}^{\|\cdot\|_{W^{3,2}}}.
\end{equation*} 
Let $\rho$ be a mollifier, i.e., $\rho\in C^\infty_c(B(0,1))$, $\rho\geq 0$, $\int_{\eR^d}\rho=1$. Then we define for $\delta>0$ 
\begin{equation}\label{MolifDef}
	\rho_\delta(\cdot)=\delta^{-d}\rho\left(\frac{\cdot}{\delta}\right).
\end{equation}
Let us notice that all the analysis presented in the paper requires the elastic energy density being of the quadratic form $W(F)=CF\cdot F+b$ for a certain fourth order tensor $C$ and $b\in\eR$. We make a minor simplification and use the elastic energy density of the form 
\begin{equation}\label{W1}
    W(F)=\frac{1}{2}|F|^2
\end{equation}
that just allows for clearer expressions but does not affect the analytical result.

We continue with the introduction of the notion of a dissipative solution. 
\begin{Definition}
Let
\begin{equation}\label{InitCondRegularity}
	u_0\in \LND(\Omega),\ F_0\in \CrlH(\Omega),\ M_0\in W^{1,2}(\Omega)^3,\ |M_0|=1\text{ a.e.\ in }\Omega.
\end{equation}
A dissipative solution to \eqref{ClassForm} with initial conditions \eqref{InitCondRegularity} and dissipation defect $\mathcal{D}\geq 0$, $\mathcal{D}\in L^\infty(0,T)$ is a triple $(u,F, M)$ enjoying the regularity
\begin{align*}
	u&\in C_w([0,T];\LND(\Omega))\cap L^2(0,T;\WND(\Omega)^d),\\
	F&\in C_w([0,T];\CrlH(\Omega)),\\
	M&\in C_w([0,T];W^{1,2}(\Omega)^3)\cap W^{1,\frac{d+2}{d+1}}(Q_T)^3
\end{align*} 
and satisfying the energy inequality
\begin{equation}\label{DSEI}
\begin{split}
	&\frac{1}{2}\left(\|u(t)\|^2_{L^2(\Omega)}+\|F(t)\|^2_{L^2(\Omega)}+\|\nabla M(t)\|^2_{L^2(\Omega)}\right)+\mathcal{D}(t)\\
	&+\int_0^t\left(\|\nabla u\|^2_{L^2(\Omega)}+\frac{1}{2}\|\tder M+(u\cdot\nabla)M\|^2_{L^2(\Omega)}\right)\\
	&\leq\int_0^t(\Hext,\tder M)+\frac{1}{2}\left(\|u_0\|^2_{L^2(\Omega)}+\|F_0\|^2_{L^2(\Omega)}+\|\nabla M_0\|^2_{L^2(\Omega)}\right)
	\end{split}
\end{equation}
for a.a.\ $t\in(0,T)$, $M$ satisfies the constraint $|M|=1$ a.e.\ in $Q_T$ and 
\begin{equation}\label{WeakForms}
\begin{split}
	(u(t),\psi(t))-(u_0,\psi(0))=&\int_0^t\Biggl((u,\tder\psi) -((u\cdot \nabla) u, \psi)-(\nabla u+FF^\top-\nabla M\odot\nabla M,\nabla\psi)\\&+\bigl((\nabla\Hext)^\top M,\psi\bigr)
	+\left\langle\mathcal{R},\nabla \psi\right\rangle\Biggr),\\
	(F(t),\Psi(t))-(F_0,\Psi(0))=&\int_0^t\Biggl((F,\tder\Psi)+(u\otimes F,\nabla \Psi)+(\nabla uF,\Psi)\Biggr),\\
	(M(t),\xi(t))-(M_0,\xi(0))=&\int_0^t \Biggl((M,\tder\xi)+\Bigl(M\times \bigl(\tder M+(u\cdot\nabla)M+2\Hext\bigr)-(u\cdot\nabla M),\xi\Bigr)\\
	&-2\sum_{i=1}^d \bigl(\partial_iM\times M,\partial_i\xi\bigr)\Biggr)
\end{split}
\end{equation}
for all $t\in(0,T)$,  $\psi\in C^1_c([0,T]\times\Omega)^d$  with $\dvr\psi=0$ in $Q_T$, $\Psi\in C^1(\overline{Q_T})^{d\times d}$, $\xi\in W^{1,2}(Q_T)^3$.
The corrector $\mathcal{R}\in L^\infty(0,T;\Meas(\overline\Omega)^{d\times d})$
satisfies
\begin{equation}\label{DefMesEst}
	\int_0^t\|\mathcal{R}(s)\|_{\Meas(\overline\Omega)}\ds\leq c\int_0^t\mathcal{D}(s)\ds.
\end{equation} 
The initial conditions are attained in the following sense
\begin{equation}\label{InitCondAtt}
\lim_{t\to 0_+}\|u(t)-u_0\|_{L^2(\Omega)}+\|F(t)-F_0\|_{L^2(\Omega)}+\|M(t)-M_0\|_{W^{1,2}(\Omega)}=0.
\end{equation}
\end{Definition}

 We remark that the integral formulation of the Landau-Lifshitz-Gilbert equation \eqref{WeakForms}$_3$ originates from a form of the Landau-Lifshitz-Gilbert equation that is equivalent to \eqref{ClassForm}$_3$, see Lemma~\ref{Lem:LLGEquivForms}.

Having all ingredients introduced we can state the main results of the paper.
\begin{Theorem}\label{Thm:Exis}
	For an arbitrary $T\in(0,\infty)$, a bounded Lipschitz domain $\Omega\subset\eR^d$, $d=2,3$, $u_0$, $F_0$ and $M_0$ satisfying \eqref{InitCondRegularity} and 
	\begin{equation}\label{HExtReg}
		\Hext\in L^2(0,T;W^{1,2}(\Omega)^3)\cap L^{d+2}(Q_T)^3
	\end{equation}
	there exists a dissipative solution to problem \eqref{ClassForm}.
\end{Theorem}
We sketch the strategy of the proof of the theorem. We begin with showing the existence of a weak solution to an approximative problem that is equipped by the original initial data for the velocity and magnetization and a regularized initial condition for the deformation gradient, cf.\ Lemma~\ref{Lem:ApproxEx}. This problem possesses an additional term regularizing the velocity in the momentum equation and the Landau-Lifshitz-Gilbert equation is replaced by a parabolic equation for an approximation $M^\varepsilon$ with a penalizing term $\varepsilon^{-1}M^\varepsilon(|M^\varepsilon|^2-1)$ that later allows to show the fulfillment of the constraint $|M|=1$. 

The existence of a solution to the approximative system is proven via a two-level Galerkin scheme. The Galerkin system contains a parabolic regularization of the transport equation for the deformation gradient yielding the compactness of approximations of the deformation gradient. Moreover, there is a cut-off function in the term involving the cross product in the equation for an approximation of the magnetization helping to bypass the lack of the constraint on the modulus of the approximation of the magnetization. 

Having the existence of approximations proven, we let the Galerkin index of the deformation gradient and the magnetization tend to infinity while keeping the Galerkin index for an approximation of the velocity fixed. At this level it is possible to adopt a proof from \cite{AlSo} to show that the modulus of approximation of the magnetization is bounded by $1$ and to remove the cut-off function from the equation. After having performed the limit passage with the Galerkin index in the approximation of the velocity,  the parabolic regularization of the transport equation vanishes but the regularity of the velocity still allows for concluding the compactness of a sequence approximating the deformation gradient. Having a solution to the approximative problem with all uniform bounds, we perform the final limit with a regularizing parameter to $0$. This then yields the existence of a dissipative solution.

\section{Proof of the main theorem}\label{Sec:DisExist}
\subsection{Approximative system}
In this section we introduce and analyze an approximating system to \eqref{ClassForm}. For $\varepsilon>0$ we consider the system
\begin{equation}\label{ApproxSysClass}
	\begin{alignedat}{2}
		\tder u+ (u\cdot \nabla) u-\Delta u-\varepsilon\Delta^3u+\nabla p-(\nabla \Hext)^\top M&=\dvr(FF^\top-\nabla M\odot\nabla M),\\
		\dvr u&=0,\\
		\tder F+\dvr(u\otimes F)&=\nabla uF,\\
		\tder M+u\cdot\nabla M-M\times\bigl(\tder M+u\cdot\nabla M+M\times\Hext\bigr)&=2\Delta M-\varepsilon^{-1}(|M|^2-1)M
	\end{alignedat}
\end{equation}
in $Q_T$ with \colorstart $\dvr F=0$\color{black}, boundary conditions 
\begin{equation*}
u=\partial_i u=\partial_i\partial_j u=0\text{ for each }i,j=1,\ldots,d,\  (\nabla M)n=0\text{ on }(0,T)\times\partial\Omega 
\end{equation*}
and initial conditions $u(0,\cdot)=u_0$, $F(0,\cdot)=F_0$, $M(0,\cdot)=M_0$ in $\Omega$ with $\dvr u_0=0$, $\dvr F_0=0$, $|M_0|\equiv 1$ in $\Omega$. The $\varepsilon$-term in the momentum equation \eqref{ApproxSysClass}$_1$ is added for regularizing the velocity, which allows for dealing with transport equation \eqref{ApproxSysClass}$_3$ in the pointwise sense, cf., Lemma~\ref{Lem:TranspEqRegData}.
We note that an approximative equation \eqref{ApproxSysClass}$_5$ with a penalizing term $M(|M|^2-1)$ is considered instead of the Landau-Lifshitz-Gilbert equation for the magnetization. The reason for this is the requirement of the low regularity of $M$, namely only square integrable gradient, that excludes the possibility to work directly with Landau-Lifshitz-Gilbert equation without additional assumptions on certain smallness of the initial data. The strategy of choosing rather a parabolic equation with the penalizing term adopts ideas from the proof of the existence result of Alouges and Soyeur in \cite{AlSo}.

The following lemma asserts the existence of a solution to system \eqref{ApproxSysClass}, which is the first task of the proof. 
\begin{Lemma}\label{Lem:ApproxEx}
	Let $\Omega\subset\eR^d$, $d=2,3$, be a bounded Lipschitz domain, $T>0$, $u_0\in L^2_{\dvr}(\Omega)^d$, $F_0\in W^{2,\infty}(\Omega)^{d\times d}$ with $\dvr F_0=0$ a.e.\ in $\Omega$, $M_0\in W^{1,2}(\Omega)^3$ with $|M_0|=1$ a.e. in $\Omega$, $\Hext$ fulfills \eqref{HExtReg}. Then there exists a weak solution to the approximative problem \eqref{ApproxSysClass}, i.e., a triple $(u, F, M)$ possessing the regularity 
	\begin{equation*}
	\begin{split}
		&u\in L^\infty (0,T;\LND(\Omega))\cap L^2(0,T;\VS(\Omega)),\ \tder u\in L^2(0,T;(\VS(\Omega))^*)\\
		&F\in C([0,T];W^{2,\infty}(\Omega)^{d\times d}),\ \tder F\in L^1(0;T;W^{1,\infty} (\Omega)^{d\times d}),\\
		& M\in L^\infty(0,T;W^{1,2}(\Omega)^3)\cap W^{1,2}(Q_T)^3
	\end{split}
	\end{equation*}
	and satisfying
	\begin{equation}\label{UFAppSys}
	\begin{alignedat}{2}
		\left\langle\tder u,\omega\right\rangle+\left((u\cdot\nabla) u,\omega\right)+\varepsilon(\nabla^3 u,\nabla^3\omega)+&\bigl(\nabla u,&&\nabla\omega\bigr)  +\bigl(FF^\top-\nabla M\odot\nabla M,\nabla\omega\bigr)\\-(\nabla\Hext^\top M,\omega)&=0 &&\text{ for all }\omega\in \VS(\Omega)\text{ and a.e.\ in }(0,T),\\
		\tder F+\dvr(u\otimes F)-\nabla uF&=0 &&\text{ a.e.\ in }Q_T,\\
		\biggl(\tder M+u\cdot\nabla M-M\times\bigl(\tder M+u\cdot\nabla M-&2M\times&&\Hext\bigr),\xi\biggr)+2(\nabla M,\nabla\xi)\\
		+\varepsilon^{-1}\bigl((|M|^2-1)M,\xi\bigr)&=0 &&\text{ for all }\xi\in W^{1,2}(\Omega)^3\text{ and a.e.\ in }(0,T).
	\end{alignedat}
	\end{equation}
	Moreover, $(u, F, M)$ fulfills
	\begin{equation}\label{EnergyBound}
	\begin{split}
		\frac{1}{2}&\left(\|u(t)\|^2_{L^2(\Omega)}+\|F(t)\|^2_{L^2(\Omega)}+\|\nabla M(t)\|^2_{L^2(\Omega)}\right)+(4\varepsilon)^{-1}\||M(t)|^2-1\|^2_{L^2(\Omega)}\\
		&+\int_0^t\|\nabla u\|^2_{L^2(\Omega)}+\varepsilon\|\nabla^3 u\|^2_{L^2(\Omega)}+\frac{1}{2}\|\tder M+(u\cdot\nabla)M\|^2_{L^2(\Omega)}\\
		\leq&\int_0^t\Bigl(\bigl(-M\times(M\times\Hext),\tder M+(u\cdot\nabla)M\bigr)+((\nabla\Hext)^\top M,u)\Bigr)\\
		&+\frac{1}{2}\left(\|u_0\|_{L^2(\Omega)}^2+\|F_0\|_{L^2(\Omega)}^2+\|\nabla M_0\|_{L^2(\Omega)}^2\right)
	\end{split}
\end{equation}
for a.a.\ $t\in(0,T)$, the constraints $\|M\|_{L^\infty(Q_T)}\leq 1$, \colorstart $\dvr F=0$  a.e.\ in $Q_T$, \color{black} and initial conditions in the sense $F(0)=F_0$, $M(0)=M_0$ and
\begin{equation}\label{InCondAttPApp}
	\lim_{t\to 0_+}\|u(t)-u_0\|_{L^2(\Omega)}=0.
\end{equation}
\end{Lemma}

 The existence proof is performed in several steps. In the first step we introduce Galerkin approximations for the velocity, the deformation gradient and the magnetization as well as an approximative system and show their existence. The system involves vanishing viscosity regularization of the equation for the deformation gradient and a parameter dependent cut-off function in the parabolic equation for $M$. We notice that adding this cut-off function is necessitated by the lack of the constraint $\|M\|_{L^\infty(Q_T)}\leq 1$ on the Galerkin level.

In the second step we collect estimates that are uniform with respect to the Galerkin index of the approximations of the deformation gradient and the magnetization; then the limit passage as this parameter tends to infinity is performed. The second step also includes the procedure originally used in \cite{AlSo} for showing that the modulus of the limit of a sequence approximating the magnetization is bounded by $1$. 

The third step consists in collecting estimates that are uniform with respect to the Galerkin index of an approximation of the velocity. We follow the procedure from \cite[Section 5.3]{BMPS18} during the derivation of the energy inequality that allows to treat the term involving the Laplacian of the second level approximation of the magnetization despite this approximation lacks square integrable second derivatives up to the boundary, which excludes the possibility of integrating by parts directly in this term. Finally, the limit passage is performed with the parameter standing for the Galerkin index for velocity approximations tending to infinity. During this limit passage, the viscosity term vanishes from the equation for the deformation gradient.

\begin{proof}
We begin with the introduction of Galerkin approximations and an approximating system. Let $\{\omega^i\}_{i=1}^\infty$ be a basis of $\VS(\Omega)$ and simultaneously an orthonormal basis of $\LND(\Omega)$. The elements of such a basis can be found as eigenfunctions to the following problem
\begin{equation*}
(\nabla^3\omega^i,\nabla^3 v)=\sigma^i(\omega^i,v)\text{ for all }v\in \VS(\Omega),
\end{equation*}
cf.\ \cite[Section A.4]{MaNeRoRu96}. We denote the projection of $\LND(\Omega)$ on $\spn\{\omega^1,\ldots,\omega^n\}$ by $P^n$. Further, let $\{\Phi^j\}_{j=1}^\infty$ and $\{\xi^n\}_{n=1}^\infty$ be orthonormal bases of $L^2(\Omega)^{d\times d}$ and $L^2(\Omega)^3$ that are orthogonal in $W^{1,2}(\Omega)^{d\times d}$ and $W^{1,2}(\Omega)^3$ respectively. We note that the existence of such bases follows from a version of the Hilbert-Schmidt theorem, cf.\ \cite[Lemma 5.1]{FiSu}. According to it the above bases consist of eigenfunctions for the following problems:
\begin{align*}
    (\nabla\Phi^j,\nabla\Theta)+(\Phi^j,\Theta)&=\Lambda^j(\Phi^j,\Theta)\text{ for all }\Theta\in W^{1,2}(\Omega)^{d\times d},\\
	(\nabla\xi^j,\nabla z)+(\xi^j,z)&=\lambda^j(\xi^j,z)\text{ for all }z\in W^{1,2}(\Omega)^3.
\end{align*}  
 The projections on finite dimensional subspaces $\spn\{\Phi^1,\ldots,\Phi^n\}$ and $\spn\{\xi^1\dots,\xi^n\}$ are denoted as $\overline{P^n}$ and $\widetilde{P^n}$, respectively. We observe that for each $n\in\eN$ it holds $\|\overline{P^n}(\Phi)\|_{W^{1,2}(\Omega)}\leq \|\Phi\|_{W^{1,2}(\Omega)}$ for all $\Phi\in W^{1,2}(\Omega)^{d\times d}$, $\|\widetilde{P^n}(\xi)\|_{W^{1,2}(\Omega)}\leq \|\xi\|_{W^{1,2}(\Omega)}$ for all $\xi\in W^{1,2}(\Omega)^3$.
 
\textbf{Step 1:} To construct the Galerkin approximations, we first introduce the following cut-off functions
\begin{align*}
	\Theta(s)&=
		\begin{cases}
			1&s\in [0,1),\\
			2-s&s\in[1,2),\\
			0&s\in [2,\infty),
		\end{cases}\\
		\Theta_k(s)&=\Theta(k^{-1}s),\ s\in[0,\infty), k\in\eN.
\end{align*} 
For fixed $m,n\in\eN$ we look for a triple $(u^{m,n},F^{m,n},M^{m,n})$ defined as
\begin{align*}
	u^{m,n}(t,x)&=\sum_{k=1}^mc^{m,n}_k(t)\omega^k(x),\ F^{m,n}(t,x)=\sum_{j=1}^nd^{m,n}_j(t)\Phi^j(x),\ M^{m,n}(t,x)=\sum_{j=1}^ne^{m,n}_j(t)\xi^j(x),
\end{align*}
where the functions $c^{m,n}=(c^{m,n}_1,\ldots, c^{m,n}_m)$, $d^{m,n}=(d^{m,n}_1,\ldots, d^{m,n}_n)$ and $e^{m,n}=(e^{m,n}_1,\ldots, e^{m,n}_n)$ satisfy, in $(0,T)$ for each $i=1,\ldots, m$ and $j=1,\ldots,n$,
\begin{equation}\label{GalSys}
	\begin{split}
		(\tder u^{m,n},\omega^i)+\bigl((u^{m,n}\cdot\nabla)u^{m,n},\omega^i\bigr)+(\nabla u^{m,n},\nabla\omega^i)+\varepsilon(\nabla^3 u^{m,n},\nabla^3 \omega^i)\\
		-\bigl(\nabla M^{m,n}\odot\nabla M^{m,n}
		-F^{m,n}(F^{m,n})^\top,\nabla\omega^i\bigr)-(\nabla\Hext^\top M^{m,n},\omega^i)&=0,\\
		(\tder F^{m,n},\Phi^j)+\bigl((u^{m,n}\cdot\nabla)F^{m,n},\Phi^j\big)-\nabla u^{m,n}F^{m,n},\Phi^j)+m^{-1}(\nabla F^{m,n},\nabla\Phi^j)&=0,\\
		(\tder M^{m,n},\xi^j)+\bigl((u^{m,n}\cdot\nabla)M^{m,n},\xi^j\bigr)-\bigl(\Theta_{m}(|M^{m,n}|)M^{m,n}\times(\tder M^{m,n}+(u^{m,n}\cdot\nabla)M^{m,n}\\
		-2M^{m,n}\times\Hext),\xi^j\bigr)
		+2(\nabla M^{m,n},\nabla\xi^j)+\varepsilon^{-1}\bigl((|M^{m,n}|^2-1)M^{m,n},\xi^j\bigr)&=0
	\end{split}
\end{equation}
with initial conditions $u^{m,n}(0)=P^m(u_0)$, $F^{m,n}(0)=\overline{P^n}(F_0)$ and $M^{m,n}(0)=\widetilde{P^n}(M_0)$. Considering now equation \eqref{GalSys}$_3$ separately and omitting the superscripts related to the indices of the approximating sequences for the rest of this step (e.g., writing $c$ instead of $c^{m,n}$), we can rewrite it in the form
\begin{equation}\label{GalMFirst}
	\tder e-A(e)\tder e=g(c,e),\ e(0)=\widetilde{P^n}M_0,
\end{equation}
with the matrix-valued function $A$ and the vector-valued function $g$ being defined as
\begin{align*}
	A_{k,l}(e)=&\sum_{r=1}^ne_r(\xi^r\times\xi^l,\xi^k),\\
	g_k(c,e)=&-\bigl((u\cdot\nabla)M,\xi^k\bigr)+\Bigl(\Theta_{m}(|M|)M\times\bigl((u\cdot\nabla)M-2M\times\Hext\bigr),\xi^k\Bigr)\\
	&-2(\nabla M,\nabla\xi^k)-\varepsilon^{-1}\bigl((|M|^2-1)M,\xi^k\bigr).
\end{align*}
We observe that $A$ is skew-symmetric due to the properties of the cross product. Hence $Id-A(e)$ is always invertible and \eqref{GalMFirst} is equivalent to 
\begin{equation}\label{GalMSecond}
	\tder e=(Id-A(e))^{-1}g(c,e),\ e(0)=\widetilde{P^n}M_0.
\end{equation}
Obviously, \eqref{GalSys}$_{1,2}$ and \eqref{GalMSecond} can be rewritten in the form of a system of $n+2m$ equations
\begin{equation*}
	\tder b=\mathcal{G}(b),\ b(0)=b_0,
\end{equation*}
where 
\begin{equation*}
	b_k=
	\begin{cases}
		c_k& k=1,\ldots,m,\\
		d_k& k=m+1,\ldots,m+n,\\
		e_k& k=m+n+1,\ldots,m+2n\\
	\end{cases}
\end{equation*}
and 
\begin{equation*}
	\mathcal{G}(b)_k=
	\begin{cases}
		-\bigl((u\cdot\nabla)u,\omega^k\bigr)-(\nabla u,\nabla\omega^k)-\varepsilon(\nabla^3 u,\nabla^3\omega^k)\\+\bigl(\nabla M\odot\nabla M-FF^\top,\nabla\omega^k\bigr)+(\nabla\Hext^\top M,\omega^k)& k=1,\ldots,m,\\
		-\bigl((u\cdot\nabla)F,\Phi^k\big)+(\nabla uF,\Phi^k)-m^{-1}(\nabla F^{m,n},\nabla\Phi^k)& k=m+1,\ldots,m+n,\\
		(Id-A(e))^{-1}g(c,e)_k& k=m+n+1,\ldots,m+2n.\\
	\end{cases}
\end{equation*}
We apply the Carath\'eodory existence theory to deduce the existence of an absolutely continuous $b$ on $(0,t^*)$ for some $t^*\in(0,T]$.

\textbf{Step 2:} Our task is the limit passage $n\to\infty$ in \eqref{GalSys} to obtain the following system for $u^m=\sum_{i=1}^mc^m_i\omega^i$, where the functions $\{\omega^i\}_{i=1}^m$ were introduced at the beginning of the proof, $F^m\in L^\infty(0,T;L^2(\Omega)^{d\times d})\cap L^2(0,T;W^{1,2}(\Omega)^{d\times d})$, $M^m\in L^2(0,T;W^{2,2}(\Omega)^3)\cap L^\infty(0,T;W^{1,2}(\Omega)^3)$ with $\tder F^m\in L^2(Q_T)^{d\times d}$ and $\tder M^m\in L^2(Q_T)^3$, 
\begin{equation}\label{MSystem}
	\begin{aligned}
		(\tder u^m,\omega^i)+\bigl((u^m\cdot\nabla)u^m,\omega^i\bigr)+(\nabla u^m,\nabla\omega^i)+\varepsilon(\nabla^3 u^m,\nabla^3\omega^i)
		&\\-\bigl(\nabla M^m\odot\nabla M^m-F^m(F^m)^\top,\nabla\omega^i\bigr)-(\nabla\Hext^\top M^m,\omega^i)&=0\\
		(\tder F^{m},\Phi)+\bigl((u^{m}\cdot\nabla)F^{m},\Phi\big)-(\nabla u^{m}F^{m},\Phi)+m^{-1}(\nabla F^{m},\nabla\Phi)&=0,\\
		\bigl((\tder M^m+(u^m\cdot\nabla)M^m-M^m\times(\tder M^m+(u^m\cdot\nabla)M^m-2M^m\times\Hext),\xi\bigr)\\
		+2(\nabla M^m,\nabla\xi)+\varepsilon^{-1}\bigl((|M^m|^2-1)M^m,\xi\bigr)&=0,
	\end{aligned}
\end{equation}
where the first identity holds a.e.\ in $(0,T)$ for all $i=1,\ldots,m,$, the second one a.e.\ in $(0,T)$ for all $\Phi\in W^{1,2}(\Omega)^{d\times d}$ and the third one a.e.\ in $(0,T)$ for all $\xi\in W^{1,2}(\Omega)^3$. Moreover, \eqref{MSystem} is accompanied with the initial conditions $u^m(0)=P^m(u_0)$, $M^m(0)=M_0$, $\lim_{t\to 0_+}\|F^m(t)-F_0\|_{L^2(\Omega)}=0$. Additionally, we have the constraint $\|M^m\|_{L^\infty(Q_T)}\leq 1$. 

We derive uniform estimates with respect to the Galerkin index $n$ but still possibly depending on $m$. First, we test \eqref{GalSys}$_3$ by $M^{m,n}$ and obtain
\begin{equation*}
\frac{1}{2}\frac{\dd}{\dt}\|M^{m,n}\|^2_{L^2(\Omega)}+2\|\nabla M^{m,n}\|^2_{L^2(\Omega)}+\varepsilon^{-1}\|M^{m,n}\|^4_{L^4(\Omega)}=\varepsilon^{-1}\|M^{m,n}\|^2_{L^2(\Omega)}.
\end{equation*}
Consequently, the application of the Gronwall lemma yields
\begin{equation}\label{MMNEst}
	\|M^{m,n}\|^2_{L^\infty(0,t^*;L^2(\Omega))}+\|\nabla M^{m,n}\|^2_{L^2(Q_{t^*})}+\|M^{m,n}\|^4_{L^4(Q_{t^*})}\leq c(\varepsilon,\|M_0\|_{L^2(\Omega)},T).
\end{equation}
Next we test \eqref{GalSys}$_2$ by $F^{m,n}$ to infer
\begin{align*}
	\frac{\dd}{\dt}\frac{1}{2}\|F^{m,n}\|^2_{L^2(\Omega)}+&\frac{1}{2}\bigl(u^{m,n},\nabla |F^{m,n}|^2\bigr)-(\nabla u^{m,n}F^{m,n},F^{m,n})\\
	+&m^{-1}\|\nabla F^{m,n}\|^2_{L^2(\Omega)}=0.
\end{align*}
Using the solenoidality of $u^{m,n}$ and the fact that $u^{m,n}=0$ on $\partial\Omega$, we get after obvious manipulations
\begin{equation}\label{FMNIneq}
	\frac{\dd}{\dt}\frac{1}{2}\|F^{m,n}\|^2_{L^2(\Omega)}+m^{-1}\|\nabla F^{m,n}\|^2_{L^2(\Omega)}-(\nabla u^{m,n}F^{m,n},F^{m,n})= 0.
\end{equation}
Testing \eqref{GalSys}$_1$ by $u^{m,n}$, we obtain 
\begin{align*}
&\frac{1}{2}\frac{\dd}{\dt}\|u^{m,n}\|_{L^2(\Omega)}^2+\|\nabla u^{m,n}\|_{L^2(\Omega)}^2+\varepsilon\|\nabla^3 u^{m,n}\|^2_{L^2(\Omega)}\\&=\bigl(\nabla M^{m,n}\odot\nabla M^{m,n}-F^{m,n}(F^{m,n})^\top,\nabla u^{m,n}\bigr)-(\nabla\Hext^\top M^{m,n},u^{m,n}).
\end{align*}
Summing up the latter identity and \eqref{FMNIneq}, we arrive at
\begin{align*}
&\frac{1}{2}\frac{\dd}{\dt}\left(\|u^{m,n}\|_{L^2(\Omega)}^2+\|F^{m,n}\|^2_{L^2(\Omega)}\right)+\|\nabla u^{m,n}\|_{L^2(\Omega)}^2+\varepsilon\|\nabla^3 u^{m,n}\|^2_{L^2(\Omega)}+ m^{-1}\|\nabla F^{m,n}\|^2_{L^2(\Omega)}\\&=(\nabla M^{m,n}\odot\nabla M^{m,n},\nabla u^{m,n})-(\nabla\Hext^\top M^{m,n},u^{m,n}).
\end{align*}
Taking into account the fact that $u^{m,n}(t)$ is an element of an $m$-dimensional space, we have
$\|u^{m,n}(t)\|_{W^{1,\infty}(\Omega)}\leq c(m)\|u^{m,n}(t)\|_{L^2(\Omega)}$. Using this fact we deduce 
\begin{align*}
\frac{1}{2}\frac{\dd}{\dt}&\bigl(\|u^{m,n}\|_{L^2(\Omega)}^2+\|F^{m,n}\|^2_{L^2(\Omega)}\bigr)+\|\nabla u^{m,n}\|_{L^2(\Omega)}^2+\varepsilon\|\nabla^3 u^{m,n}\|^2_{L^2(\Omega)}+ m^{-1}\|\nabla F^{m,n}\|^2_{L^2(\Omega)}\\
\leq &\|\nabla M^{m,n}\|^2_{L^2(\Omega)}\|u^{m,n}\|_{W^{1,\infty}(\Omega)}+\|\nabla\Hext\|_{L^2(\Omega)} \|M^{m,n}\|_{L^2(\Omega)}\|u^{m,n}\|_{L^\infty(\Omega)}\\
\leq &c(m)\left(\|\nabla M^{m,n}\|^2_{L^2(\Omega)}+\|\nabla\Hext\|_{L^2(\Omega)} \|M^{m,n}\|_{L^2(\Omega)}\right)\|u^{m,n}\|_{L^2(\Omega)}\\
\leq &c(m)\left(\|\nabla\Hext\|^2_{L^2(\Omega)}\|M^{m,n}\|^2_{L^2(\Omega)}+\|\nabla M^{m,n}\|^2_{L^2(\Omega)}+\left(1+\|\nabla M^{m,n}\|^2_{L^2(\Omega)}\right)\|u^{m,n}\|^2_{L^2(\Omega)}\right).
\end{align*}
Using \eqref{MMNEst}, we obtain by the Gronwall lemma
\begin{equation}\label{UFMNEst}
	\begin{split}
	&\|u^{m,n}\|^2_{L^\infty(0,t^*;L^2(\Omega))}+\|F^{m,n}\|^2_{L^\infty(0,t^*;L^2(\Omega))}+\|u^{m,n}\|^2_{L^2(0,t^*;W^{1,2}(\Omega))}+ m^{-1}\|\nabla F^{m,n}\|^2_{L^2(Q_{t^*})}\\
	&\leq c\left(m,\varepsilon,\|u_0\|_{L^2(\Omega)},\|F_0\|_{L^2(\Omega)},\|M_0\|_{L^2(\Omega)},\|\Hext\|_{L^2(0,T;W^{1,2}(\Omega))},T\right).
	\end{split}
\end{equation}
Hence we deduce that $t^*=T$ and 
\begin{equation}\label{UMNEst}
	\|u^{m,n}\|_{L^\infty(0,T;W^{1,\infty}(\Omega))}\leq c(m)\|u^{m,n}\|_{L^\infty(0,t^*;L^2(\Omega))}=c(m)\|c^{m,n}\|_{L^\infty(0,t^*)}\leq c.
\end{equation}
The next task is to derive uniform estimates on time derivatives of $F^{m,n}$, $M^{m,n}$ and $u^{m,n}$. 
Testing \eqref{GalSys}$_2$ by $\tder F^{m,n}$, we obtain
\begin{equation*}
	\|\tder F^{m,n}\|^2_{L^2(\Omega)}+(2m)^{-1}\frac{\dd}{\dt}\|\nabla F^{m,n}\|^2_{L^2(\Omega)}=-\bigl((u^{m,n}\cdot\nabla )F^{m,n},\tder F^{m,n}\bigr)+(\nabla u^{m,n}F^{m,n},\tder F^{m,n}).
\end{equation*}
An estimate of the right hand side of the latter equality by the Young inequality yields
\begin{equation*}
\|\tder F^{m,n}\|^2_{L^2(\Omega)}+(2m)^{-1}\frac{\dd}{\dt}\|\nabla F^{m,n}\|^2_{L^2(\Omega)}\leq c\|u^{m,n}\|^2_{W^{1,\infty}(\Omega)}\|F^{m,n}\|^2_{W^{1,2}(\Omega)}.
\end{equation*}
Integrating over $(0,t)\subset (0,T)$ and using \eqref{UFMNEst}, we have
\begin{equation}\label{TDerFMNEst}
	\|\tder F^{m,n}\|^2_{L^2(Q_T)}+\|\nabla F^{m,n}\|^2_{L^\infty(0,T;L^2(\Omega))}\leq c(m,\varepsilon,u_0,F_0,M_0,\Hext,T). 
\end{equation}
Testing \eqref{GalSys}$_3$ by $\tder M^{m,n}$, we obtain
\begin{equation}\label{TDerMMNTestId}
	\begin{split}
		&\|\tder M^{m,n}\|^2_{L^2(\Omega)}+\frac{\dd}{\dt}\left(\|\nabla M^{m,n}\|^2_{L^2(\Omega)}+(4\varepsilon)^{-1}\||M^{m,n}|^2-1\|^2_{L^2(\Omega)}\right)\\
		&=\Bigl(-(u^{m,n}\cdot\nabla)M^{m,n}+\Theta_{m}(|M^{m,n}|)M^{m,n}\times\bigl((u^{m,n}\cdot\nabla)M^{m,n}-2M^{m,n}\times\Hext\bigr),\tder M^{m,n}\Bigr).
	\end{split}
\end{equation}
Employing the Young inequality on the right hand side of the latter identity, we get
\begin{align*}
	\|\tder M^{m,n}\|^2_{L^2(\Omega)}+\frac{\dd}{\dt}\bigl(\|\nabla M^{m,n}\|^2_{L^2(\Omega)}&+(4\varepsilon)^{-1}\||M^{m,n}|^2-1\|^2_{L^2(\Omega)}\bigr)\\
	&\leq c\|u^{m,n}\|^2_{L^\infty(\Omega)}\|\nabla M^{m,n}\|^2_{L^2(\Omega)}(1+m^2)+cm^2\|\Hext\|^2_{L^2(\Omega)}.
\end{align*}
Applying the Gronwall lemma, \eqref{MMNEst} and \eqref{UMNEst}, we arrive at
\begin{equation}\label{TDerMMNEst}
\begin{split}
	\|\tder M^{m,n}\|^2_{L^2(Q_T)}&+\|\nabla M^{m,n}\|^2_{L^\infty(0,T;L^2(\Omega))}+(2\varepsilon)^{-1}\||M^{m,n}|^2-1\|^2_{L^\infty(0,T;L^2(\Omega))}\\&\leq c(m,\varepsilon,u_0,F_0,M_0,\Hext,T).
	\end{split}
\end{equation}
We note that in order to estimate the term $\||M^{m,n}(0)|^2-1\|^2_{L^2(\Omega)}$ appearing during the above computations, one applies the embedding $W^{1,2}(\Omega)$ to $L^4(\Omega)$ and the continuity of the projection $\widetilde{P^n}$ to get
\begin{equation*}
\||M^{m,n}(0)|^2-1\|^2_{L^2(\Omega)}\leq c(1+\|M^{m,n}(0)\|^4_{L^{4}(\Omega)})\leq c(1+\|M^{m,n}(0)\|^4_{W^{1,2}(\Omega)})\leq c(1+\|M_0\|^4_{W^{1,2}(\Omega)}).
\end{equation*}
It follows directly from \eqref{GalSys}$_1$ thanks to \eqref{MMNEst}, \eqref{UFMNEst} and \eqref{TDerMMNEst} that
\begin{equation}\label{CMNTDerEst}
	\|\tder c^{m,n}\|_{L^\infty(0,T)}\leq c(m,\varepsilon,u_0,F_0,M_0,\Hext,T).
\end{equation}
 We are ready to derive convergences that are essential for the passage $n\to\infty$. As a direct consequence of estimates \eqref{UMNEst}, \eqref{CMNTDerEst} and the Arzel\`a-Ascoli  and Banach-Alaoglu theorems one obtains the existence of a not explicitly labeled subsequence such that 
\begin{equation}\label{CMNConv}
	\begin{alignedat}{2}
		c^{m,n}&\rightharpoonup^* c^m &&\text{ in }W^{1,\infty}(0,T;\eR^m),\\
		c^{m,n}&\to c^m&&\text{ in }C([0,T];\eR^m).
	\end{alignedat}
\end{equation}
Hence by the definition of $u^{m,n}$ we immediately get 
\begin{equation}\label{UMNConv}
\begin{alignedat}{2}
	u^{m,n}&\rightharpoonup^* u^m &&\text{ in }W^{1,\infty}(0,T;W^{1,\infty}(\Omega)^d),\\
	u^{m,n}&\to u^m&&\text{ in }C([0,T];W^{1,\infty}(\Omega)^d),
\end{alignedat}
\end{equation}
where
\begin{equation*}
	u^m(t,x)=\sum_{i=1}^m c_i^{m}(t)\omega^i(x),\ u^{m}(0,x)=P^m(u_0)(x).
\end{equation*}
 Thanks to \eqref{TDerFMNEst}, \eqref{TDerMMNEst} and the Aubin-Lions lemma along with \eqref{UFMNEst} and \eqref{MMNEst}, we deduce the existence of $(c^m,F^m,M^m)$ and not explicitly labeled subsequences such that
\begin{equation}\label{FMMNConv}
	\begin{alignedat}{2}
		\tder F^{m,n}&\rightharpoonup \tder F^m&&\text{ in }L^2(Q_T)^{d\times d},\\
		F^{m,n}&\to F^m&&\text{ in }L^2(0,T;L^2(\Omega)^{d\times d}),\\
		F^{m,n}&\rightharpoonup F^m&&\text{ in }L^2(0,T;W^{1,2}(\Omega)^{d\times d}),\\
		M^{m,n}&\rightharpoonup^*M^m&&\text{ in }L^\infty(0,T; W^{1,2}(\Omega)^3),\\
		\tder M^{m,n}&\rightharpoonup \tder M^m&&\text{ in }L^2(Q_T)^3,\\
		M^{m,n}&\to M^m&&\text{ in }L^2(0,T;L^q(\Omega)^3),\ q\in[1,2^*)\text{ and a.e.\ in }Q_T.
	\end{alignedat}
\end{equation}
Using the fact that $\Theta_m(|M^{m,n}|)=0$ for $|M^{m,n}|>2m$, the convergence $M^{m,n}\to M^m$ a.e.\ in $Q_T$ and the Lebesgue dominated convergence theorem, we get
\begin{equation}\label{CutMNConv}
	\Theta_m(|M^{m,n}|)M^{m,n}\to \Theta_m(|M^m|)M^m \text{ in }L^r(Q_T),\ r\in [1,\infty).
\end{equation}
Next we observe that by the interpolation with the bound on $M^{m,n}$ in $L^\infty(0,T;L^2(\Omega)^3)$ we obtain from convergence \eqref{FMMNConv}$_6$
\begin{equation}\label{MNL4Conv}
	M^{m,n}\to M^m\text{ in }L^4(Q_T)^3.
\end{equation}
Now we are ready to pass to the limit $n\to\infty$ in \eqref{GalSys}. Let us fix $m\in\eN$, $\psi\in C^1_c(0,T)$ and $\xi^j$. Then we multiply \eqref{GalSys}$_3$ by $\psi$, integrate the result over $(0,T)$ and perform the passage $n\to\infty$ with the help of \eqref{FMMNConv}$_{4,5,6}$, \eqref{UMNConv}$_2$,  \eqref{CutMNConv} and \eqref{MNL4Conv} to eventually arrive at
\begin{equation}\label{MMIdent}
	\begin{split}
&(\tder M^m,\xi)+\bigl((u^m\cdot\nabla)M^m,\xi\bigr)-\bigl(\Theta_{m}(|M^m|)M^m\times(\tder M^m+(u^m\cdot\nabla)M^m-2M^m\times\Hext),\xi\bigr)\\
		&+2(\nabla M^m,\nabla\xi)+\varepsilon^{-1}(|M^m|^2-1)M^m,\xi\bigr)=0\text{ a.e.\ in }(0,T)
		\end{split}
\end{equation}
for arbitrary $\xi \in W^{1,2}(\Omega)^3$ since $\{\xi^j\}_{j=1}^\infty$ forms a basis in $W^{1,2}(\Omega)^3$. 

The next task concerning this equation is to show that $\|M^m\|_{L^\infty(Q_T)}\leq 1$, which allows us to remove the cut-off function $\Theta_m(|M^m|)$. Let the functions $g,G:\eR\to\eR$ be defined as 
\begin{equation*}
g(s)=
\begin{cases}
	0&s\in(-\infty,0),\\
	s&s\in[0,1),\\
	1&s\in[1,\infty),
\end{cases} 
\quad
G(s)=
\begin{cases}
	0&s\in(-\infty,0),\\
	\frac{s^2}{2}&s\in[0,1),\\
	s-\frac{1}{2}&s\in[1,\infty).
\end{cases} 
\end{equation*}
Note that $G$ is a primitive function to $g$. We set $\eta:=g(|M^m|^2-1)M^m$. Then we observe that $\eta\in W^{1,2}(Q_T)^3$ as $g(|M^m|^2-1)$ $g'(|M^m|^2-1)|M^m|^2\in L^\infty(Q_T)$ due to $g'(|M^m|^2-1)=0$ if $|M^m|^2>2$. Hence setting $\xi:=\eta(t)$ in \eqref{MMIdent} and integrating the result over $(0,t)$ for fixed $t$, we obtain
\begin{equation}\label{LLGgTest}
	\begin{split}
&\int_{Q_t}\bigl(\tder M^m+(u^m\cdot\nabla)M^m\bigr)\cdot g(|M^m|^2-1)M^m
		+\int_{Q_t}2\partial_i M^m\cdot\bigl(g(|M^m|^2-1)\partial_i M^m\\
		&+2g'(|M^m|^2-1)(M^m\cdot\partial_i M^m)M^m\bigr)+\varepsilon^{-1}(|M^m|^2-1)|M^m|^2g(|M^m|^2-1)=0.
		\end{split}
\end{equation}
We rewrite the first integral in \eqref{LLGgTest} with the help of the function $G$, the solenoidality of $u^m$ and $u^m=0$ on $\partial\Omega$ in the following way:
\begin{align*}
\int_{Q_t}\bigl(\tder M^m+(u^m\cdot\nabla)M^m\bigr)\cdot g(|M^m|^2-1)M^m&=\frac{1}{2}\int_{Q_t}\tder G(|M^m|^2-1)+u^m\cdot\nabla G(|M^m|^2-1)\\
&=\frac{1}{2}\int_{Q_t}\tder G(|M^m|^2-1).
\end{align*}
The next observation is that the integrand of the second integral in \eqref{LLGgTest} is nonnegative. Hence we infer
\begin{equation*}
	\int_{Q_t}\tder G(|M^m|^2-1)\leq 0,
\end{equation*}
which implies that for an arbitrary $t\in(0,T)$ we have
\begin{equation*}
	\int_{\Omega}G(|M^m(t)|^2-1)\leq \int_{\Omega}G(|M^m(0)|^2-1)=\int_{\Omega}G(|M_0|^2-1)=0.
\end{equation*}
Hence we conclude $|M^m|\leq 1$ a.e.\ in $Q_T$. It follows that $\Theta_m(|M^m|)=1$ a.e.\ in $Q_T$ and combining this fact with \eqref{MMIdent} we obtain \eqref{MSystem}$_3$. We note that setting $\xi:=M^m(t)$ in \eqref{MSystem}$_3$ and integrating over $(0,T)$ one arrives at 
\begin{equation}\label{MmEqTest}
\int_{Q_T} |\nabla M^m|^2=-\frac{1}{2}\int_0^T \Bigl((\tder M^m,M^m)+\varepsilon^{-1}\bigl((|M^m|^2-1)M^m,M^m\bigr)\Bigr).
\end{equation}
Our intention now is to show 
\begin{equation}\label{NablaMNStrong}
    \nabla M^{m,n}\to \nabla M^m\text{ in }L^2(Q_T)^{3\times d}.
\end{equation}
To this end we employ \eqref{GalSys}$_3$ to infer
\begin{equation*}
    \lim_{n\to\infty}\int_{Q_T}|\nabla M^{m,n}|^2=-\frac{1}{2}\int_0^T\Bigl((\tder M^m,M^m)+\varepsilon^{-1}\bigl((|M^m|^2-1)M^m,M^m\bigr)\Bigr)=\int_{Q_T} |\nabla M^m|^2
\end{equation*}
by \eqref{FMMNConv}$_5$, \eqref{MNL4Conv} and \eqref{MmEqTest}. To conclude \eqref{NablaMNStrong} it suffices to combine the latter convergence of norms with \eqref{FMMNConv}$_4$. 
Next we we fix $j\in\eN$ and $\psi\in C^\infty_c(0,T)$, multiply \eqref{GalSys}$_2$ by $\psi$, integrate  the identity over $(0,T)$ and pass to the limit $n\to\infty$ with the help of \eqref{UMNConv}$_2$ and \eqref{FMMNConv}$_{1,2,3}$ to get
\begin{equation}\label{FMIdent}
	\begin{split}
	\int_0^T\psi\bigl((\tder F^m+(u^m\cdot\nabla)F^m+\nabla u^mF^m,\Phi^j)+m^{-1}(\nabla F^m,\nabla\Phi^j)\bigr)=0.
	\end{split}
\end{equation}
As $\{\Phi^j\}_{j=1}^\infty$ forms a basis in $W^{1,2}(\Omega)^{d\times d}$, we can replace $\Phi^j$ in \eqref{FMIdent} by an arbitrary function from $W^{1,2}(\Omega)^{d\times d}$. We conclude \eqref{MSystem}$_2$ by the fundamental theorem of the calculus of variations.

The last task is the limit passage $n\to\infty$ in the equality \eqref{GalSys}$_1$. We multiply \eqref{GalSys}$_1$ for fixed $i\in\{1,\ldots,m\}$ by an arbitrary but fixed $\psi\in C^\infty_c(0,T)$, integrate the resulting identity over $(0,T)$ and pass to the limit $n\to\infty$ therein using convergences \eqref{CMNConv}, \eqref{FMMNConv}$_2$, \eqref{MNL4Conv} and \eqref{NablaMNStrong} to obtain
\begin{align*}
	&\int_0^T\psi \bigl((\tder u^m,\omega^i)+\bigl((u^m\cdot\nabla)u^m,\omega^i\bigr)+(\nabla u^{m},\nabla\omega^i)+(\nabla^3 u^{m},\nabla^3\omega^i)\\&
		-\bigl(\nabla M^{m}\odot\nabla M^{m}-F^{m}(F^{m})^\top,\nabla\omega^i\bigr)-(\nabla\Hext^\top M^{m},\omega^i)\bigr)=0.
\end{align*}
We note that due to \eqref{FMMNConv}$_2$ we get $F^{m,n}(F^{m,n})^\top\to F^m(F^m)^\top$ in $L^1(Q_T)^{d\times d}$ as $n\to\infty$. We conclude \eqref{MSystem}$_1$ by the fundamental theorem of calculus of variations. The attainment of the initial conditions by $F^m$, $M^m$ follows directly as $\tder F^m\in L^2(0,T;L^2(\Omega)^{d\times d})$ and the construction yields $F^m(0)=F_0$, $M^m(0)=M_0$ a.e.\ in $\Omega$.

\textbf{Step 3}: We pass to the limit $m\to\infty$ in \eqref{MSystem} to obtain \eqref{UFAppSys} for the limit functions $(u,F,M)$ possessing the regularity expressed in the statement of the lemma. We begin with collecting estimates that are independent of $m$ and essential for deducing necessary convergences for the limit passage $m\to \infty$. We would like to multiply \eqref{MIdent2} by $\tder M^m+(u^m\cdot\nabla)M^m$ and integrate the resulting identity over $\Omega$ to get further estimates and terms that will be later canceled after testing \eqref{MSystem}$_1$ by $u^m$. However, the available regularity of $M^m$ excludes the possibility of performing the integration by parts in the term $\int_\Omega\Delta M^m\cdot\tder M^m$. To circumvent this inconvenience we integrate \eqref{TDerMMNTestId} over $(0,t)\subset(0,T)$ to get
\begin{align*}
		&\|\nabla M^{m,n}(t)\|^2_{L^2(\Omega)}+(4\varepsilon)^{-1}\||M^{m,n}(t)|^2-1\|^2_{L^2(\Omega)}+\int_0^t\|\tder M^{m,n}\|^2_{L^2(\Omega)}\\
		&=\int_0^t\Biggl(\Bigl(-(u^{m,n}\cdot\nabla)M^{m,n}+\Theta_{m}(|M^{m,n}|)M^{m,n}\times\bigl((u^{m,n}\cdot\nabla)M^{m,n}-2M^{m,n}\times\Hext\bigr),\tder M^{m,n}\Bigr)\Biggr)\\
		&\quad+\|\nabla M^{m,n}(0)\|^2_{L^2(\Omega)}+(4\varepsilon)^{-1}\||M^{m,n}(0)|^2-1\|^2_{L^2(\Omega)}.
	\end{align*}
Then using \eqref{HExtReg}, \eqref{UMNConv}$_2$, \eqref{FMMNConv}$_{5,7}$, \eqref{MNL4Conv}, \eqref{NablaMNStrong}, the weak lower semicontinuity of norms, the convergence of $M^{m,n}(0)$ towards $M_0$ in $L^4(\Omega)^3$, which is a consequence of the embedding $W^{1,2}(\Omega)^3$ to $L^4(\Omega)^3$ and the fact that $M^{m,n}(0)\to M_0$ in $W^{1,2}(\Omega)^3$ as $n\to\infty$ and $|M_0|=1$ a.e.\ in $\Omega$, we conclude
\begin{equation}\label{TDerMMTest}
	\begin{split}
			&\|\nabla M^m(t)\|^2_{L^2(\Omega)}+(4\varepsilon)^{-1}\||M^m(t)|^2-1\|^2_{L^2(\Omega)}+\int_0^t\|\tder M^m\|^2_{L^2(\Omega)}\\
		&\leq \int_0^t\Bigl(-(u^m\cdot\nabla)M^m+M^m\times\bigl((u^m\cdot\nabla)M^m-2M^m\times\Hext\bigr),\tder M^m\Bigr)+\|\nabla M_0\|^2_{L^2(\Omega)}
	\end{split}
\end{equation}
first for a.a.\ $t\in(0,T)$ and later for all $t\in(0,T)$ as $M^m$ is in fact continuous in time with respect to the weak topology of $W^{1,2}(\Omega)^3$. We note that the strong convergence 
\begin{equation*}
\Theta_{m}(|M^{m,n}|)M^{m,n}\times(u^{m,n}\cdot\nabla)M^{m,n}\to\Theta_{m}(|M^{m}|)M^{m}\times(u^{m}\cdot\nabla)M^{m}\text{ in }L^2(Q_t)^3
\end{equation*}
follows by the generalized dominated convergence theorem as $|\Theta_{m}(|M^{m,n}|)M^{m,n}\times(u^{m,n}\cdot\nabla)M^{m,n}|\leq 2m|u^{m,n}||\nabla M^{m,n}|$ and the term on the right hand side of the latter inequality converges towards $2m|u^{m}||\nabla M^{m}|$ in $L^2(Q_t)$ due to \eqref{UMNConv}$_2$ and \eqref{NablaMNStrong}. Next, due to the regularity of the limit functions in \eqref{FMMNConv}$_{4,5}$ and \eqref{UMNConv}, we have $\Delta M^{m}\in L^2(Q_T)^3$ from \eqref{MMIdent} and by the interior elliptic regularity
\begin{equation}\label{MMInterReg}
    M^m\in L^2(0,T;W^{2,2}(\tilde\Omega)^3)\text{ for each }\tilde\Omega\subset\subset\Omega.
\end{equation}
The choice of an arbitrary $\xi\in C^\infty_c(\Omega)^3$ in \eqref{MMIdent} and the fundamental theorem of the calculus of variations then yield
\begin{equation}\label{MIdent2}
\begin{split}
&\tder M^m+(u^m\cdot\nabla)M^m-\Theta_{m}(|M^m|)M^m\times\bigl(\tder M^m+(u^m\cdot\nabla)M^m-2M^m\times\Hext\bigr)\\
&-2\Delta M^m+\varepsilon^{-1}(|M^m|^2-1)M^m=0\text{ a.e.\ in }Q_T.
\end{split}
\end{equation}
We multiply \eqref{MIdent2} by $(u^m\cdot\nabla)M^m$ and integrate over $\Omega$ to get
\begin{equation}\label{ConvDerMTest}
\begin{split}
	(\tder M^m,(u^m\cdot\nabla)M^m)+\|(u^m\cdot\nabla)M^m\|^2_{L^2(\Omega)}&=\bigl(M^m\times(\tder M^m-2M^m\times\Hext),(u^m\cdot\nabla)M^m\bigr)\\&+\bigl(2\Delta M^m-\varepsilon^{-1}(|M^m|^2-1)M^m,(u^m\cdot\nabla)M^m\bigr).
	\end{split}
\end{equation}
We handle the last term on the right hand of the latter equality. Splitting it into two parts 
\begin{equation*}
I_1=2\bigl((u^m\cdot\nabla)M^m,\Delta M^m\bigr),\ I_2=-\varepsilon^{-1}\bigl((u^m\cdot\nabla)M^m,M^m(|M^m|^2-1)\bigr)
\end{equation*}
we have immediately
\begin{equation*}
I_2=-(4\varepsilon)^{-1}\int_\Omega u^m\cdot\nabla(|M^m|^2-1)^2=0
\end{equation*}
as $u^m$ is divergence free and its trace vanishes in $(0,T)\times\partial\Omega$.
As the available regularity of $M^m$ does not guarantee the integrability of the second derivatives up to $\partial\Omega$, cf., \eqref{MMInterReg}, one has to handle $I_1$ more carefully than just to integrate by parts. Considering a smooth function $\kappa_\delta$, $0\leq\kappa_\delta\leq 1$ such that
\begin{equation*}
    \kappa_\delta(x)=
    \begin{cases}
        0&\dist(x,\partial\Omega)\leq \delta,\\
        1&\dist(x,\partial\Omega)\geq 2\delta
    \end{cases}
\end{equation*}
and $|\nabla\kappa_\delta|\leq c\delta^{-1}$, we can write 
\begin{equation}\label{DeltaLim}
    I_1=2\int_{\Omega }(u^m\cdot\nabla)M^m\Delta M^m=\lim_{\delta\to 0_+}2\int_{\Omega }(u^m\cdot\nabla)M^m\Delta M^m\kappa_\delta,
\end{equation}
since the left integral is well defined. For the integral under the limit we further compute, employing the Einstein convention,
\begin{align*}
    \int_{\Omega }&(u^m)_i\partial_i(M^m)_k\partial_{ll}(M^m)_k\kappa_\delta\\
    =&-\int_{\Omega }\partial_l(u^m)_i\partial_i(M^m)_k\partial_l(M^m)_k\kappa_\delta+(u^m)_i\partial_{il}(M^m)_k\partial_l(M^m)_k\kappa_\delta+(u^m)_i\partial_i(M^m)_k\partial_l(M^m)_k\partial_l\kappa_\delta\\
    =&-\int_{\Omega }\partial_l(u^m)_i\partial_i(M^m)_k\partial_l(M^m)_k\kappa_\delta+\frac{1}{2}(u^m)_i\partial_{i}\partial_l(M^m)^2_k\kappa_\delta+(u^m)_i\partial_i(M^m)_k\partial_l(M^m)_k\partial_l\kappa_\delta\\
    =&-\int_{\Omega }\partial_l(u^m)_i\partial_i(M^m)_k\partial_l(M^m)_k\kappa_\delta-\frac{1}{2}(u^m)_i\partial_l(M^m)^2_k\partial_i\kappa_\delta+(u^m)_i\partial_i(M^m)_k\partial_l(M^m)_k\partial_l\kappa_\delta\\
    =&-\int_{\Omega }\partial_i(M^m)_k\partial_l(M^m)_k \partial_l(u^m)_i\kappa_\delta\\
    &+\int_{\{0<\kappa_\delta<1\}}\frac{1}{2}(u^m)_i[\partial_l(M^m)_k]^2\partial_i\kappa_\delta-(u^m)_i\partial_i(M^m)_k\partial_l(M^m)_k\partial_l\kappa_\delta
\end{align*}
using \eqref{MMInterReg}, the integration by parts, the fact that $\kappa_\delta$ has a compact support in $\Omega$ and the solenoidality of $u^m$. The latter computation gives
\begin{equation}\label{ProdIneq}
\begin{split}
    &\int_{\Omega }(u^m\cdot\nabla )M^m\Delta M^m\kappa_\delta\\
    &=-\int_{\Omega }\nabla M^m\odot M^m\cdot \nabla u^m\kappa_\delta+\int_{\{0<\kappa_\delta<1\}}\frac{1}{2}|\nabla M^m|^2u^m\cdot\nabla \kappa_\delta-(\nabla M^m\odot\nabla M^m)\cdot(u^m\otimes \nabla\kappa_\delta).
    \end{split}
\end{equation}
In order to proceed, we estimate for $x\in\Omega$ with $\dist(x,\partial\Omega)\leq 2\delta$
\begin{equation*}
    |u^m(x)||\nabla\kappa_\delta(x)|=|u^m(x)-u^m(\tilde x)||\nabla\kappa_\delta(x)|\leq c \|\nabla u^m\|_{L^\infty(\Omega)}|x-\tilde x|\delta^{-1},
\end{equation*}
where $\tilde x\in \partial\Omega$ is such that $|x-\tilde x|=\dist(x,\partial\Omega)$, $u^m(\tilde x)=0$ accordingly. We note that, as $\Omega$ is Lipschitz, we can always find such a point $\tilde x$. Hence we obtain
\begin{equation}\label{DeltaEst}
    \sup_{\{\kappa_\delta\in (0,1)\}}|u^m||\nabla\kappa_\delta|\leq c,
\end{equation}
where the constant is independent of $\delta$. Recalling that $\nabla M^m(t)\in L^2(\Omega)^{3\times d}$ for a.e.\  $t\in(0,T)$, we combine \eqref{ProdIneq} with \eqref{DeltaEst} to get 
\begin{equation*}
    I_1=\lim_{\delta\to 0_+}2\int_{\Omega }(u^m\cdot\nabla)M^m\Delta M^m\kappa_\delta=\lim_{\delta\to 0_+}-2\int_{\Omega }\nabla M^m\odot M^m\cdot \nabla u^m\kappa_\delta=-2\int_{\Omega }\nabla M^m\odot M^m\cdot \nabla u^m
\end{equation*}
from \eqref{DeltaLim}.
Integrating \eqref{ConvDerMTest} over $(0,t)$ using the results of the above computations and adding \eqref{TDerMMTest}, we deduce
\begin{equation}\label{ConvDerMMTest}
\begin{split}
	&\|\nabla M^m(t)\|^2_{L^2(\Omega)}+(4\varepsilon)^{-1}\||M^m(t)|^2-1\|^2_{L^2(\Omega)}+\int_0^t\|\tder M^m+(u^m\cdot\nabla)M^m\|^2_{L^2(\Omega)}\\
	&\leq 2\int_0^t\Bigl(-\bigl(M^m\times(M^m\times\Hext),\tder M^m+(u^m\cdot\nabla)M^m\bigr)-(\nabla M^m\odot\nabla M^m,\nabla u^m)\Bigr)+\|\nabla M_0\|^2_{L^2(\Omega)}.
\end{split}	
\end{equation}
Testing \eqref{MSystem}$_1$ by $u^m$ and \eqref{MSystem}$_2$ by $F^m$, we obtain similarly as in the step 2
\begin{equation}\label{UMTest}
\begin{split}
	&\frac{1}{2}\frac{\dd}{\dt}\left(\|u^m\|^2_{L^2(\Omega)}+\|F^m\|^2_{L^2(\Omega)}\right)+\|\nabla u^m\|^2_{L^2(\Omega)}+ m^{-1}\|\nabla F^m\|^2_{L^2(\Omega)}+\varepsilon\|\nabla^3 u^m\|^2_{L^2(\Omega)}\\
	&=(\nabla M^m\odot\nabla M^m,\nabla u^m)+((\nabla\Hext)^\top M^m,u^m).
\end{split}
\end{equation}
Integrating \eqref{UMTest} over $(0,t)$, adding the result to \eqref{ConvDerMMTest} multiplied by $\frac{1}{2}$, we arrive at
\begin{equation}\label{EnergyIneqM}
\begin{split}
	\frac{1}{2}&\left(\|u^m(t)\|^2_{L^2(\Omega)}+\|F^m(t)\|^2_{L^2(\Omega)}+\|\nabla M^m(t)\|^2_{L^2(\Omega)}+(4\varepsilon)^{-1}\||M^m(t)|^2-1\|^2_{L^2(\Omega)}\right)\\
	&+\int_0^t\bigl(\|\nabla u^m\|^2_{L^2(\Omega)}+\varepsilon\|\nabla^3 u^m\|^2_{L^2(\Omega)}+ m^{-1}\|\nabla F^m\|^2_{L^2(\Omega)}+\frac{1}{2}\|\tder M^m+(u^m\cdot\nabla)M^m\|^2_{L^2(\Omega)}\bigr)\\
	\leq &\int_0^t\Bigl(\bigl(-M^m\times(M^m\times\Hext),\tder M^m+(u^m\cdot\nabla)M^m\bigr)+((\nabla\Hext)^\top M^m,u^m)\Bigr)\\&
	+\frac{1}{2}\left(\|u^m(0)\|^2_{L^2(\Omega)}+\|F_0\|^2_{L^2(\Omega)}+\|\nabla M_0\|^2_{L^2(\Omega)}\right).
\end{split}
\end{equation}
Using $\|M^m\|_{L^\infty(Q_T)}\leq 1$ and the Poincar\'e and Young inequalities, we infer 
\begin{equation}\label{EpsEnergEst}
\begin{split}
	\|&u^m(t)\|^2_{L^2(\Omega)}+\|F^m(t)\|^2_{L^(\Omega)}+\|\nabla M^m(t)\|^2_{L^2(\Omega)}+(2\varepsilon)^{-1}\||M^m(t)|^2-1\|^2_{L^2(\Omega)}\\
	&+\int_0^t\left(2\left(\|\nabla u^m\|^2_{L^2(\Omega)}+\varepsilon\|\nabla^3 u^m\|^2_{L^2(\Omega)}+m^{-1}\|\nabla F^m\|^2_{L^2(\Omega)}\right)+\|\tder M^m+(u^m\cdot\nabla)M^m\|^2_{L^2(\Omega)}\right)\\
	\leq& c\bigl(1+\|\Hext\|^2_{L^2(0,T;W^{1,2}(\Omega))}+\|u_0\|^2_{L^2(\Omega)}+\|F_0\|^2_{L^2(\Omega)}+\|\nabla M_0\|^2_{L^2(\Omega)}\bigr).
\end{split}
\end{equation}
The latter inequality implies the following bounds
\begin{equation}\label{MEst}
\begin{split}
	&\|u^m\|^2_{L^\infty(0,T;L^2(\Omega))}+\|u^m\|^2_{L^2(0,T;W^{3,2}(\Omega))}+m^{-1}\|\nabla F^m\|^2_{L^2(Q_T)}+\|M^m\|^2_{L^\infty(0,T;W^{1,2}(\Omega))}\\&+\|\tder M^m+(u^m\cdot\nabla)M^m\|^2_{L^2(Q_T)}\leq c(\varepsilon,u_0,F_0,M_0,\Hext,T).
	\end{split}
\end{equation}
We note that the bound on the $W^{3,2}$-norm of $u^m$ follows as $\|\nabla^3\cdot  \|_{L^2(\Omega)}$ is equivalent to the standard $W^{3,2}$-norm on $\VS(\Omega)$ by the Friedrichs inequality. 
For the time derivative $\tder M^m$ we then have due to the embedding $W^{3,2}(\Omega)$ to $L^\infty(\Omega)$
\begin{equation}\label{TDerMMEst}
\begin{split}
\|\tder M^m\|_{L^2(Q_T)}\leq &\|\tder M^m+(u^m\cdot\nabla)M^m\|_{L^2(Q_T)}+\||u^m||\nabla M^m|\|_{L^2(Q_T)}\\
\leq &c+\|u^m\|_{L^2(0,T;L^\infty(\Omega))}\|M^m\|_{L^\infty(0,T;W^{1,2}(\Omega))}\leq c(\varepsilon,u_0,F_0,M_0,\Hext,T).
\end{split}
\end{equation}
It follows from \eqref{MSystem}$_1$ by using \eqref{MEst} and $\|M^m\|_{L^\infty(Q_T)}\leq 1$ that
\begin{equation}\label{TderUMEst}
	\|\tder u^m\|_{L^2(0,T;(\VS(\Omega))^*)}\leq c(\varepsilon,u_0,F_0,M_0,\Hext,T).
\end{equation}
Invoking standard compactness arguments, cf.\ Step 2, we deduce from \eqref{MEst}, \eqref{TDerMMEst} and \eqref{TderUMEst} the existence of a triple $(u,F,M)$ and not explicitly labeled subsequences such that
\begin{equation}\label{MConv}
	\begin{alignedat}{2}
		u^m&\rightharpoonup^* u&&\text{ in }L^\infty(0,T;L^2(\Omega)^d),\\
		u^m&\rightharpoonup u&&\text{ in }L^2(0,T;W^{3,2}(\Omega)^d),\\
		u^m&\to u&&\text{ in }L^2(0,T;W^{1,\infty}(\Omega)^d)\text{ and a.e.\ in }Q_T,\\
		\tder u^m&\rightharpoonup \tder u&&\text{ in }L^2(0,T;(\mathcal{V}(\Omega))^*),\\
		F^m&\rightharpoonup^* F&&\text{ in }L^\infty(0,T;L^2(\Omega)^{d\times d}),\\
		M^m&\rightharpoonup^* M&&\text{ in }L^\infty(0,T;W^{1,2}(\Omega)^d),\\
		M^m&\to M&&\text{ in }L^2(0,T;L^q(\Omega)^3),\ q\in[1,2^*),\\
		\tder M^m&\rightharpoonup \tder M&&\text{ in }L^2(Q_T)^3,\\
		M^m&\to M&&\text{ in }L^p(Q_T)^3,\ p\in[1,\infty),\\
		\tder M^m+(u^m\cdot\nabla)M^m&\rightharpoonup \tder M+(u\cdot\nabla)M &&\text{ in }L^2(Q_T)^3.
	\end{alignedat}
\end{equation}
We note that \eqref{MConv}$_{9}$ is a direct consequence of $\|M^m\|_{L^\infty(Q_T)}\leq 1$ and  \eqref{MConv}$_{7}$. Next, \eqref{MConv}$_{10}$ follows from \eqref{MEst} and \eqref{MConv}$_{3,6}$. 

Let us begin with the limit passage $m\to\infty$ in \eqref{MSystem}. 
We fix $\xi\in W^{1,2}(\Omega)^3$, multiply \eqref{MSystem}$_3$ by a fixed $\psi\in C^\infty_c(0,T)$, integrate the result over $(0,T)$ and pass to the limit $m\to\infty$ with the help of \eqref{MConv}$_{6,9,10}$. We arrive at
\begin{equation}\label{MEqWeakly}
\begin{split}
	\int_0^T\psi\Bigl(\tder M+(u\cdot\nabla)M&-M\times\bigl(\tder M+(u\cdot\nabla)M-2M\times\Hext\bigr),\xi\Bigr)\\
	&=-\int_0^T\psi\bigl(2(\nabla M,\nabla\xi)+\varepsilon^{-1}(|M|^2-1)M,\xi)\bigr),
	\end{split}
\end{equation}
from which \eqref{UFAppSys}$_3$ follows. Moreover, from identity \eqref{MMIdent} we obtain
\begin{equation*}
\lim_{m\to\infty}\int_{Q_T}|\nabla M^m|^2=-\frac{1}{2}\int_0^T\Bigl((\tder M,M)+\varepsilon^{-1}\bigl((|M|^2-1)M,M\bigr)\Bigr)=\int_{Q_T}|\nabla M|^2
\end{equation*}
using \eqref{MConv}$_{8,9}$ and \eqref{MEqWeakly} with $\xi:=M$. Combining the latter equalities with \eqref{MConv}$_6$ we have
\begin{equation}\label{NablaMMStrongly}
    \nabla M^m\to \nabla M\text{ in }L^2(Q_T)^{3\times d}.
\end{equation}
Next we fix $\Psi\in C^1(\overline{Q_T})^{d\times d}$, $s\in(0,T)$, multiply \eqref{MSystem}$_2$ by $\Psi(s)$, integrate over $Q_t$, $t\in(0,T)$, integrate by parts in time and space variables, multiply the result by an arbitrary $\psi\in C^\infty_c(0,T)$ and integrate over $(0,T)$ to obtain
\begin{align*}
	\int_0^T\psi(t)\bigl((F^m(t),\Psi(t))-(F_0,\Psi(0))\bigr)=&\int_0^T\psi(t)\int_0^t\bigl((F^m,\tder\Psi)+(u^m\otimes F^m,\nabla\Psi)\\
	&+(\nabla u^m F^m,\Psi)-m^{-1}(\nabla F^m,\nabla\Psi)\bigr).
\end{align*}
Hence passing to the limit $m\to\infty$ in the latter identity using \eqref{MConv}$_{3,5}$ and \eqref{MEst} we arrive at
\begin{equation*}
	(F(t),\Psi(t))-(F_0,\Psi(0))=\int_0^t\bigl((F,\tder\Psi)+(u\otimes F,\nabla\Psi)+(\nabla u F,\Psi)\bigr).
\end{equation*}
As a consequence of the amount of the regularity possessed by $u$ and $F_0$ it follows from Lemma~\ref{Lem:TranspEqRegData} that $F\in C([0,T];W^{2,\infty}(\Omega)^{d\times d})$ with $\tder F\in L^1(0,T;W^{1,\infty}(\Omega)^{d\times d}))$ and $F$ satisfies
\eqref{UFAppSys}$_2$. 

The next task is to show that in fact $F^m\to F$ in $L^\infty(0,T;L^2(\Omega)^{d\times d})$. To this end we fix $\Phi\in W^{1,2}(\Omega)^{d\times d}$, multiply \eqref{UFAppSys}$_2$ by $\Phi$, integrate over $\Omega$, subtract the result from \eqref{MSystem}$_2$, set $\Phi=F^m-F$, integrate the result over $(0,t)\subset(0,T)$ and add to both sides of the result $-m^{-1}\int_0^t\bigl(\nabla F,\nabla (F^m-F)\bigr)$ to obtain
\begin{equation}\label{FMFDiffIneq}
\begin{split}
	&\|(F^m-F)(t)\|^2_{L^2(\Omega)}+m^{-1}\|\nabla (F^m-F)\|^2_{L^2(Q_t)}=\int_0^t\Bigl(\bigl((u^m-u)\cdot\nabla\bigr)F,F^m-F\Bigr)\\&
	+\bigl(\nabla u^m(F^m-F)+(\nabla u^m-\nabla u)F,F^m-F\bigr)-m^{-1}\int_0^t\bigl(\nabla F,\nabla (F^m-F)\bigr)=\sum_{k=1}^3 I^m_k.
\end{split}
\end{equation}
We estimate separately each term on the right hand side of the latter identity. By the obvious manipulations and the Young inequality we get 
\begin{align*}
	|I^m_1|\leq & c\int_0^t\Bigl(\|u^m-u\|^2_{L^\infty(\Omega)}\|\nabla F\|^2_{L^\infty(Q_T)}+\|F^m-F\|^2_{L^2(\Omega)}\Bigr),\\
	|I^m_2|\leq  & c\int_0^t\Bigl(\|\nabla u^m\|_{L^\infty(Q_T)}\|F^m-F\|^2_{L^2(\Omega)}+\|\nabla u^m-\nabla u\|^2_{L^2(\Omega)}\|F\|^2_{L^\infty(Q_T)}+\|F^m-F\|^2_{L^2(\Omega)}\Bigr).
\end{align*}
In order to estimate $I^m_3$ we first integrate by parts. We note that at this moment we do not know that the normal derivative of $F$ on $\partial\Omega$ vanishes; therefore we also have to estimate the boundary integral and obtain
\begin{align*}
|I^m_3|\leq& m^{-1}\Biggl(\Bigl|\int_0^t(\Delta F,F^m-F)\Bigr|+\Bigl|\int_0^t\int_{\partial\Omega}(\nabla F)n\cdot (F^m-F)\Bigr|\Biggr)\\\leq& c\int_0^tm^{-2}\|\Delta F\|^2_{L^2(\Omega)}+\|F^m-F\|^2_{L^2(\Omega)}
+m^{-1}\int_0^t\|\nabla F\|_{L^2(\partial\Omega)}\|F^m-F\|_{L^2(\partial\Omega)}.
\end{align*}
Using the properties of the trace operator and the Young inequality in the last term, we deduce
\begin{align*}
m^{-1}\int_0^t\|\nabla F\|_{L^2(\partial\Omega)}\|F^m-F\|_{L^2(\partial\Omega)}\leq& \int_0^tcm^{-1}\|\nabla F\|^2_{W^{1,2}(\Omega)}\\&+(2m)^{-1}\left(\|F^m-F\|^2_{L^2(\Omega)}+\|\nabla(F^m-F)\|^2_{L^2(\Omega)}\right).
\end{align*}
Employing the above estimates in \eqref{FMFDiffIneq}, we get an inequality for $t\in(0,T)$
\begin{equation*}
	\|(F^m-F)(t)\|^2_{L^2(\Omega)}+(2m)^{-1}\|\nabla (F^m-F)\|^2_{L^2(Q_t)}\leq \int_0^th^m(s)+c\|(F^m-F)(s)\|^2_{L^2(\Omega)}\ds,
\end{equation*}
where 
\begin{equation*}
	h^m=c\bigl(m^{-1}+\|u^m-u\|^2_{L^\infty(\Omega)}+\|\nabla u^m-\nabla u\|^2_{L^2(\Omega)}\bigr)\|F\|^2_{L^\infty(0,T;W^{2,\infty}(\Omega))}
\end{equation*}
and $\int_0^T h^m(s)\ds\to 0$ as $m\to\infty$ as a consequence of \eqref{MConv}$_{3}$. Then we conclude that
\begin{equation}\label{FMStrongC}
	F^m\to F\text{ in }L^\infty(0,T;L^2(\Omega)^{d\times d})\text{ as }m\to\infty
\end{equation}
by the Gronwall lemma. Having the latter convergence at hand we deduce $F^m(F^m)^\top \to FF^\top$ in $L^\infty(0,T;L^1(\Omega)^{d\times d})$ as $m\to\infty$. Thus we are ready for the passage $m\to\infty$ in \eqref{MSystem}$_1$ multiplied by $\phi\in C^\infty_c(0,T)$ and integrated over $(0,T)$. To conclude \eqref{UFAppSys}$_1$ we also employ \eqref{MConv}$_{2,4}$, \eqref{NablaMMStrongly}, \eqref{FMStrongC} and the fact that $\{\omega^i\}_{i=1}^\infty$ is a basis in $\mathcal{V}(\Omega)$.

Finally, we focus on showing inequality \eqref{EnergyBound}. Multiplying \eqref{EnergyIneqM} by $\theta\in C^\infty_c(0,T)$, $\theta\geq 0$ and integrating over $(0,T)$, we get 
\begin{equation*}
\begin{split}
	\int_0^T&\theta(\tau)\left(\frac{1}{2}\left(\|u^m(\tau)\|^2_{L^2(\Omega)}+\|F^m(\tau)\|^2_{L^2(\Omega)}+\|\nabla M^m(\tau)\|^2_{L^2(\Omega)}\right)\right.\\&\left.+\int_0^\tau\|\nabla u^m\|^2_{L^2(\Omega)}+\varepsilon\|\nabla^3 u^m\|^2_{L^2(\Omega)}+\frac{1}{2}\|\tder M^m+(u^m\cdot\nabla)M^m\|^2_{L^2(\Omega)}\right)\dd\tau\\ 
	\leq& \int_0^T\theta(\tau)\int_0^\tau \bigl(-M^m\times(M^m\times\Hext),\tder M^m+(u^m\cdot\nabla)M^m\bigr)+(\nabla\Hext^\top M^m,u^m)\\
	&+\frac{1}{2}\left(\|u_0\|^2_{L^2(\Omega)}+\|F_0\|^2_{L^2(\Omega)}+\|\nabla M_0\|^2_{L^2(\Omega)}\right)\dd\tau.
	\end{split}
\end{equation*}
The passage $m\to\infty$ in the latter inequality with the help of convergences \eqref{MConv}$_{1,2,3,5,6,9,10}$, \eqref{HExtReg} and the weak lower semicontinuity of the norm yield
\begin{equation*}
\begin{split}
	\int_0^T&\theta(\tau)\left(\frac{1}{2}\left(\|u(\tau)\|^2_{L^2(\Omega)}+\|F(\tau)\|^2_{L^2(\Omega)}+\|\nabla M(\tau)\|^2_{L^2(\Omega)}\right)\right.\\&\left.+\int_0^\tau\|\nabla u\|^2_{L^2(\Omega)}+\varepsilon\|\nabla^3 u\|^2_{L^2(\Omega)}+\frac{1}{2}\|\tder M+(u\cdot\nabla)M\|^2_{L^2(\Omega)}\right)\dd\tau\\
	\leq& \int_0^T\theta(\tau)\int_0^\tau \bigl(-M\times(M\times\Hext),\tder M+(u\cdot\nabla)M\bigr)+(\nabla\Hext^\top M,u)\\
	&+\frac{1}{2}\left(\|u_0\|^2_{L^2(\Omega)}+\|F_0\|^2_{L^2(\Omega)}+\|\nabla M_0\|^2_{L^2(\Omega)}\right)\dd\tau.
	\end{split}
\end{equation*}
Fixing $t\in(0,T)$ and setting  $\theta(\tau)=\rho_\delta(t-\tau)$ in the latter inequality, where $\rho_\delta$ is a one-dimensional mollifier with $\delta<\frac{1}{2}\min\{t,T-t\}$, and letting $\delta\to 0_+$ we conclude \eqref{EnergyBound}.
The attainment of the initial condition $u_0$ in the sense \eqref{InCondAttPApp} follows immediately as $u\in L^2(0,T;\mathcal{V}(\Omega))$ and $\tder u\in L^2(0,T;(\mathcal{V}(\Omega))^*)$ implies $u\in C(0,T;L^2(\Omega)^d)$. Further, $u(0)=u_0$ follows in a standard way.
\end{proof}
After having proved the existence of solutions to approximative system~\eqref{ApproxSysClass} and after having shown various estimates of the solutions that are independent of the regularizing parameter $\varepsilon$, we next focus on  several convergence results.

\begin{Lemma}\label{Lem:Convergences}
	Let $\Omega$, $u_0$, $F_0$, $M_0$, $W$ and $\Hext$ satisfy the assumptions of Lemma~\ref{Lem:ApproxEx}. Let $\{\varepsilon^r\}_{r=1}^\infty$ be a sequence such that $\varepsilon^r\to 0_+$ as $r\to\infty$ and $\{(u^r, F^r, M^r)\}_{r=1}^\infty$ be a sequence of weak solutions to \eqref{ApproxSysClass} with $\varepsilon=\varepsilon^r$ constructed in Lemma~\ref{Lem:ApproxEx}. Then the following uniform estimates hold 
\begin{equation}\label{RUnifEst}
	\begin{split}
		\|u^r\|_{L^\infty(0,T;L^2(\Omega))}&\leq c,\\
		\|u^r\|_{L^2(0,T;W^{1,2}(\Omega))}&\leq c,\\
		(\varepsilon^r)^\frac{1}{2}\|\nabla^3 u^r\|_{L^2(0,T;L^2(\Omega))}&\leq c,\\
		\|F^r\|_{L^\infty(0,T;L^2(\Omega))}&\leq c,\\
		\|M^r\|_{L^\infty(0,T;W^{1,2}(\Omega))}&\leq c,\\
		(\varepsilon^r)^{-\frac{1}{2}}\||M^r|^2-1\|_{L^\infty(0,T;L^2(\Omega))}&\leq c,\\
		\|\tder M^r+(u^r\cdot\nabla)M^r\|_{L^2(Q_T)}&\leq c,\\
		\|M^r\|_{L^\infty(Q_T)}&\leq 1,\\
		\|\tder u^r\|_{L^2(0,T;(\mathcal{V}(\Omega))^*)}&\leq c,\\
		\|\tder M^r\|_{L^\frac{d+2}{d+1}(Q_T)}&\leq c,
	\end{split}
\end{equation}
and there exist a not explicitly labeled subsequence of $\{(u^r, F^r, M^r)\}_{r=1}^\infty$, $u\in L^\infty(0,T;\LND(\Omega))\cap L^2(0,T;\WND(\Omega)^d)$, $F\in L^\infty(0,T;\CrlH(\Omega))$, $M\in L^\infty(0,T;W^{1,2}(\Omega)^3)$, $\overline{|F|^2}\in L^\infty(0,T;\NnMeas(\overline{\Omega}))$ and $\overline{|G|^2}\in L^\infty(0,T;\NnMeas(\overline{\Omega}))$ such that	
\begin{equation}\label{NConvergences}
		\begin{alignedat}{2}
			u^r&\rightharpoonup u&&\text{ in }L^2(0,T;\WND(\Omega)^d),\\
			u^r&\rightharpoonup^* u&&\text{ in }L^\infty(0,T;L^2(\Omega)^d),\\
			u^r&\to u&&\text{ in }L^2(0,T;L^2(\Omega)^d),\\
			\tder u^r&\rightharpoonup \tder u&&\text{ in }L^2(0,T;(\mathcal{V}(\Omega)^*),\\
			F^r&\rightharpoonup^* F&&\text{ in }L^\infty(0,T;L^2(\Omega)^{d\times d}),\\
			\tder F^r&\rightharpoonup \tder F&&\text{ in }L^2(0,T;(W^{3,2}(\Omega)^{d\times d})^*),\\
			M^r&\rightharpoonup^* M&&\text{ in }L^\infty(0,T;W^{1,2}(\Omega)^3),\\
			M^r&\to M&&\text{ in }L^p(Q_T)^3,\ p\in[1,\infty),\\
			|M^r|^2&\to 1&&\text{ in }L^\infty(0,T;L^2(\Omega)),\\
			\tder M^r&\rightharpoonup \tder M&&\text{ in }L^\frac{d+2}{d+1}(Q_T)^3,\\
			\tder M^r+(u^r\cdot\nabla)M^r&\rightharpoonup\tder M+(u\cdot\nabla)M&&\text{ in }L^2(Q_T)^3,\\
			|F^r|^2&\rightharpoonup^* \overline{|F|^2}&&\text{ in }L^\infty(0,T;\NnMeas(\overline{\Omega})),\\
			|\nabla M^r|^2&\rightharpoonup^* \overline{|G|^2}&&\text{ in }L^\infty(0,T;\NnMeas(\overline{\Omega})).
		\end{alignedat}
	\end{equation}
	\begin{proof}
		Let $\{(u^r,F^r, M^r)\}_{r=1}^\infty$ be a sequence of solutions to \eqref{ApproxSysClass} from the assertion of the lemma. By \eqref{EnergyBound} and the constraint $\|M^r\|_{L^\infty(Q_T)}\leq 1$ we immediately obtain the estimates in \eqref{RUnifEst}$_{1,2,3,4,5,6,7}$. Moreover, \eqref{UFAppSys}$_1$ yields, for arbitrary $\phi\in L^2(0,T;\VS(\Omega))$, that 
		\begin{equation*}
		\begin{split}
			\Biggl|\int_0^T\left\langle \tder u^r,\phi\right\rangle\Biggr|
			=&\Biggl|\int_0^T(u^r\otimes u^r-\nabla u^r+\nabla M^r\odot\nabla M^r-F^r(F^r)^\top,\nabla\phi)\\
			&-\varepsilon^r(\nabla^3 u^r,\nabla^3 \phi)+\bigl((\nabla\Hext)^\top M^r,\phi\bigr)\Biggr| \\
			\leq &c\left(\|u^r\|^2_{L^\infty(0,T;L^2(\Omega))}+\|F^r\|^2_{L^\infty(0,T;L^2(\Omega))}+\|M^r\|^2_{L^\infty(0,T;L^2(\Omega))} \right.\\
			&\left. +\|\nabla u^r\|_{L^2(Q_T)}\right)\|\nabla\phi\|_{L^2(0,T;L^\infty(\Omega))}+c\varepsilon^r\|\nabla^3 u^r\|_{L^2(0,T;L^2(\Omega))}\|\phi\|_{L^2(0,T;W^{3,2}(\Omega))}\\
			&+\|\Hext\|_{L^2(0,T;W^{1,2}(\Omega))}\|M^r\|_{L^\infty(Q_T)}\|\phi\|_{L^2(Q_T)}\\
			\leq &c(\|u_0\|^2_{L^2(\Omega)}, \|F_0\|^2_{L^2(\Omega)}, \|M_0\|^2_{L^2(\Omega)}, \|\Hext\|_{L^2(0,T;W^{1,2}(\Omega))})\bigl(1+(\varepsilon^r)^\frac{1}{2}\bigr)\|\phi\|_{L^2(0,T;W^{3,2}(\Omega)}.
			\end{split}
		\end{equation*}
		Since $\{\varepsilon^r\}_{r=1}^\infty$ is bounded as a convergent sequence, we deduce 
		\begin{equation}\label{TDerUE}
		\|\tder u^r\|_{L^2(0,T;(\mathcal{V}(\Omega))^*)}\leq c.
		\end{equation}
		Fixing $\Phi\in L^2(0,T;W^{3,2}(\Omega)^{d\times d})$ we obtain from \eqref{UFAppSys}$_2$
		\begin{equation*}
		\begin{split}
			\Biggl|\int_0^T\left\langle \tder F^r,\Phi\right\rangle\Biggr|&= \Biggl|\int_0^T\int_\Omega (u^r\otimes F^r)\cdot\nabla\Phi+\nabla u^rF^r\cdot\Phi\Biggr|\\
			&\leq c\|u^r\|_{L^2(0,T;W^{1,2}(\Omega))}\|F^r\|_{L^\infty(0,T;L^2(\Omega))}\|\Phi\|_{L^2(0,T;W^{3,2}(\Omega))}.
		\end{split}
		\end{equation*}
		Hence we conclude
		\begin{equation} \label{TDerFE}
		\|\tder F^r\|_{L^2(0,T;(W^{3,2}(\Omega))^*)}\leq c.
		\end{equation}
		As a consequence of the Gagliardo-Nirenberg interpolation theorem we get
		\begin{equation*}
			\|u^r\|^\frac{2(d+2)}{d}_{L^\frac{2(d+2)}{d}(Q_T)}\leq c\|u^r\|^\frac{4}{d}_{L^\infty(0,T;L^2(\Omega))}\|u^r\|^2_{L^2(0,T;W^{1,2}(\Omega))}.
		\end{equation*}
		Accordingly, we obtain by the H\"older inequality and \eqref{RUnifEst}$_{1,2,5,7}$
		\begin{equation}\label{TDerMrEst}
			\begin{split}
			\|\tder M^r\|_{L^\frac{d+2}{d+1}(Q_T)}\leq &\|\tder M^r+(u^r\cdot\nabla)M^r\|_{L^\frac{d+2}{d+1}(Q_T)}+\||u^r||\nabla M^r|\|_{L^\frac{d+2}{d+1}(Q_T)}\\
			\leq &c \|\tder M^r+(u^r\cdot\nabla)M^r\|_{L^2(Q_T)}+\|u^r\|_{L^\frac{2(d+2)}{d}(Q_T)}\|\nabla M^r\|_{L^2(Q_T)}\leq c.
			\end{split}
		\end{equation}
The weak($^*$) convergences in \eqref{NConvergences} and convergence \eqref{NConvergences}$_9$ are obtained as a direct consequence of \eqref{RUnifEst}, \eqref{TDerUE} and \eqref{TDerMrEst} whereas \eqref{NConvergences}$_3$ and \eqref{NConvergences}$_8$ follow by the Aubin-Lions lemma. Note that in order to show the latter convergence also the constraint \eqref{RUnifEst}$_8$ is applied and \eqref{NConvergences}$_{11}$ follows due to \eqref{NConvergences}$_{3,7,10}$.
	\end{proof}
\end{Lemma}

\subsection{Proof of Theorem~\ref{Thm:Exis}}
\textbf{Step 1}: For given initial data, we first construct a sequence of weak solutions to approximative problem \eqref{ApproxSysClass}. Let $\{\varepsilon^r\}_{r=1}^\infty$ be a sequence such that $\varepsilon^r\to 0_+$ as $r\to \infty$. Next we define $F^r_0$ as the mollification of $F_0$, i.e., $(F^r_0)_{ij}=(F_0)_{ij}*\rho_{\varepsilon^r}$, $i,j=1\ldots,d$ with $\rho_{\varepsilon^r}$ as in \eqref{MolifDef}. Then we have $F^r_0\in W^{2,\infty}(\Omega)^{d\times d}$ and $\dvr F^r_0=0$ in $\Omega$. Applying Lemma~\ref{Lem:Convergences} with the initial data $u_0$, $F^r_0$ and $M_0$, we find a sequence $\{(u^r,F^r,M^r)\}_{r=1}^\infty$ of solutions to \eqref{ApproxSysClass} with $\varepsilon=\varepsilon^r$ and a limit triple $(u,F,M)$.

\textbf{Step 2}: We derive the energy inequality for $(u,F,M)$. From the energy inequality \eqref{EnergyBound} we infer for a fixed $\tau\in(0,T)$ and $r\in\eN$ 
\begin{equation*}
\begin{split}
	\frac{1}{2}\int_\Omega&\left(|u(\tau)|^2+|F(\tau)|^2+|\nabla M(\tau)|^2\right)+\int_0^\tau\int_\Omega\left(|\nabla u^r|^2+\frac{1}{2}|\tder M^r+(u^r\cdot\nabla)M^r|^2\right)\\&+\frac{1}{2}\int_\Omega\left(|\nabla M^r(\tau)|^2-|\nabla M(\tau)|^2+|F^r(\tau)|^2-|F(\tau)|^2\right)\\
	\leq& \int_0^\tau \bigl(-M^r\times(M^r\times\Hext),\tder M^r+(u^r\cdot\nabla)M^r\bigr)+(\nabla\Hext^\top M^r,u^r)\\
	&+\frac{1}{2}\int_\Omega\left(|u_0|^2+|F^r_0|^2+|\nabla M_0|^2\right).
	\end{split}
\end{equation*}
Using the identity $a\times(b\times c)=(a\cdot c)b-(a\cdot b)c$ in the latter inequality and multiplying the result by $\theta\in C^\infty_c(0,T)$, $\theta\geq 0$ and integrating over $(0,T)$, we get 
\begin{equation*}
\begin{split}
	\int_0^T&\theta(\tau)\left(\frac{1}{2}\int_\Omega\left(|u(\tau)|^2+|F(\tau)|^2+|\nabla M(\tau)|^2\right)+\int_0^\tau\int_\Omega\left(|\nabla u^r|^2+\frac{1}{2}|\tder M^r+(u^r\cdot\nabla)M^r|^2\right)\right)\dd\tau\\&+\frac{1}{2}\int_0^T\theta(\tau)\int_\Omega\left(|F^r(\tau)|^2-|F(\tau)|^2+|\nabla M^r(\tau)|^2-|\nabla M(\tau)|^2\right)\dd\tau\\ \leq& \int_0^T\theta(\tau)\left(\int_0^\tau \bigl(-(M^r\cdot\Hext)M^r+|M^r|^2\Hext,\tder M^r+(u^r\cdot\nabla)M^r\bigr)+(\nabla\Hext^\top M^r,u^r)\right.\\
	&\left.+\frac{1}{2}\int_\Omega\left(|u_0|^2+|F_0^r|^2+|\nabla M_0|^2\right)\right)\dd\tau.
	\end{split}
\end{equation*}
Before we pass to the limit $r\to\infty$, we observe that
\begin{equation*}
	\bigl(\Hext,(u\cdot\nabla)M\bigr)=\int_\Omega {(\Hext)}_iu_j\partial_j M_i=-\int_\Omega \partial_j{(\Hext)}_iu_jM_i=-(\nabla\Hext^\top M,u).
\end{equation*}
Now, using \eqref{NConvergences}, the weak lower semicontinuity of the norm and the fact that $|M|=1$ a.e.\ in $Q_T$;  we arrive, in the limit as $r\to \infty$, at
\begin{equation}\label{EILimWeak}
\begin{split}
	\int_0^T&\theta(\tau)\left(\frac{1}{2}\int_\Omega\left(|u(\tau)|^2+|F(\tau)|^2+|\nabla M(\tau)|^2\right)+\int_0^\tau\int_\Omega\left(|\nabla u|^2+\frac{1}{2}|\tder M+(u\cdot\nabla)M|^2\right)+\mathcal{D}(\tau)\right)\dd\tau\\
	\leq &\int_0^T\theta(\tau)\left(\int_0^\tau (\Hext,\tder M)+\frac{1}{2}\int_\Omega\left(|u_0|^2+|F_0|^2+|\nabla M_0|^2\right)\right)\dd\tau,
	\end{split}
\end{equation}
where the dissipation defect $\mathcal{D}$ is defined as
\begin{equation*}
	\mathcal{D}(t)= \left(\overline{|F(t)|^2}-|F(t)|^2\dx\right)(\overline\Omega)+\left(\overline{|G(t)|^2}-|\nabla M(t)|^2\dx\right)(\overline\Omega).
\end{equation*}
The regularity of the limit objects $\overline{|F|^2}$, $\overline{|G|^2}$, $F$ and $M$ yields that $\mathcal{D}\in L^\infty(0,T)$. The nonnegativity of $\mathcal{D}$ follows since convex functionals are weak lower semicontinuous. Fix $t\in(0,T)$ and set $\theta(\tau)=\rho_\delta(t-\tau)$ in \eqref{EILimWeak}, where $\rho_\delta$ is a one-dimensional mollifier with $\delta<\frac{1}{2}\min\{t,T-t\}$. In the limit as $\delta\to 0_+$ we obtain
\begin{align*}
	\frac{1}{2}&\left(\|u(t)\|_{L^2(\Omega)}^2+\|F(t)\|_{L^2(\Omega)}^2+\|\nabla M(t)\|_{L^2(\Omega)}^2\right)+\mathcal{D}(t)+\int_0^t\left(\|\nabla u\|_{L^2(\Omega)}^2+\frac{1}{2}\|\tder M+(u\cdot\nabla)M\|^2_{L^2(\Omega)}\right)\\
	\leq& \int_0^t(\Hext,\tder M)+\frac{1}{2}\left(\|u_0\|_{L^2(\Omega)}^2+\|F_0\|_{L^2(\Omega)}^2 +\|\nabla M_0\|_{L^2(\Omega)}^2\right)
\end{align*}
for a.a.\ $t\in (0,T)$. 

\textbf{Step 3}: We pass to the limit $r\to\infty$ in the formulation of the approximative problem. We begin with the convergence of the sequences $\{F^r(F^r)^\top\}_{r=1}^\infty$, $\{\nabla M^r\odot\nabla M^r\}_{r=1}^\infty$ for which only an $L^1$ uniform estimate with respect to the space variable is available. From \eqref{RUnifEst}$_4$ we infer the existence of a not explicitly labeled subsequence $\{F^r(F^r)^\top\}_{r=1}^\infty$
and $\mathcal{R}_1\in L^\infty(0,T;\Meas(\overline\Omega)^{d\times d})$ such that
\begin{equation}\label{NonlinConvNFF}
		F^r(F^r)^\top\rightharpoonup^* FF^\top+\mathcal{R}_1\text{ in }L^\infty(0,T;\Meas(\overline\Omega)^{d\times d})\text{ as }r\to\infty.
\end{equation}
We estimate the corrector $\mathcal{R}_1$ with the help of the dissipation defect $\mathcal{D}$.
Fix $\Phi\in C(\overline{\Omega})^{d\times d}$ arbitrarily. We then get for a.a.\ $t\in(0,T)$
\begin{equation*}
	\begin{split}
	\int_0^t &\left\langle F^r(F^r)^\top-FF^\top,\Phi\right\rangle\ds\\=&\int_0^t \left\langle(F^r-F)(F^r-F)^\top,\Phi\right\rangle +\left\langle F(F^r-F)^\top,\Phi\right\rangle+\left\langle(F^r-F)F^\top,\Phi\right\rangle\ds\\ \leq& \int_0^t \left\langle|F^r-F|^2,|\Phi|\right\rangle+\int_\Omega \left((F^r-F)\cdot\Phi^\top F+(F^r-F)\cdot\Phi F\right)\ds.
	\end{split}
\end{equation*}
Next we let $r$ tend to $\infty$ and employ \eqref{NConvergences}$_{5,12}$ and \eqref{NonlinConvNFF} to obtain
\begin{equation}\label{FirstDef}
	\int_0^t \left\langle \mathcal{R}_1,\Phi\right\rangle\ds\leq \int_0^t \left\langle \overline{|F|^2}-|F|^2,|\Phi|\right\rangle\ds.
	\end{equation}
Similarly we get 
\begin{equation}\label{NonlinConvNMNM}
		\nabla M^r\odot\nabla M^r\rightharpoonup^* \nabla M\odot\nabla M+\mathcal{R}_2\text{ in }L^\infty(0,T;\Meas(\overline\Omega)^{d\times d})\text{ as }r\to\infty
\end{equation}
and for a.a.\ $t\in(0,T)$ by \eqref{NConvergences}$_{7,13}$
\begin{equation}\label{SecDef}
	\begin{split}
	\int_0^t \left\langle \mathcal{R}_2,\Phi\right\rangle=&
	\lim_{r\to\infty}\int_0^t \left\langle \nabla M^r\odot\nabla M^r-\nabla M\odot\nabla M,\Phi\right\rangle\ds\\
	=&\lim_{r\to\infty}\int_0^t \left\langle\nabla (M^r-M)\odot\nabla(M^r-M),\Phi\right\rangle+\left\langle\nabla M\odot\nabla (M^r-M)\right.\\
	&\left.+\nabla (M^r-M)\odot\nabla M,\Phi\right\rangle\ds\\
	\leq &\limsup_{r\to\infty}\int_0^t \left\langle|\nabla (M^r-M)|^2,|\Phi|\right\rangle=\int_0^t \left\langle\overline{|G|^2}-|\nabla M|^2,|\Phi|\right\rangle.
	\end{split}
\end{equation}
Next we take the supremum over $\Phi\in C(\overline{\Omega})^{d\times d}$ with $\|\Phi\|_{C(\overline\Omega)}\leq 1$ in \eqref{FirstDef} and \eqref{SecDef} and deduce that 
\begin{equation*}
		\int_0^t \|\mathcal{R}_1+\mathcal{R}_2\|_{\Meas(\overline{\Omega})}\ds\leq c\int_0^t\mathcal{D}(s)\ds\text{ for a.a.\ }t\in(0,T).
\end{equation*}

We are now in a position to show that the integral formulations in \eqref{WeakForms} hold true. We first consider \eqref{WeakForms}$_1$. Fixing $s\in(0,T)$ and setting $\omega=\psi(s)$ in \eqref{ApproxSysClass}$_1$, where $\psi\in C^1_c([0,T]\times\Omega)^d$, $\dvr\psi=0$ in $Q_T$, integrating the result over $(0,t)$, integrating by parts in time, applying convergences \eqref{NConvergences}$_{1,3,4}$, \eqref{NonlinConvNFF} and \eqref{NonlinConvNMNM}, employing estimate \eqref{RUnifEst}$_3$ and setting $\mathcal{R}=\mathcal{R}_1+\mathcal{R}_2$, we conclude \eqref{WeakForms}$_1$.

We focus on the passage $\varepsilon\to 0$ in \eqref{ApproxSysClass}$_2$ leading to \eqref{WeakForms}$_2$. Multiplying \eqref{ApproxSysClass}$_2$ by $\Phi\in C^1(\overline{Q_T})^{d\times d}$ and integrating over $Q_t$, $t\in(0,T)$, we obtain after an integration by parts
\begin{equation}\label{RTransport}
	(F^r(t),\Phi(t))-(F^r_0,\Phi(0))=\int_0^t(F^r,\tder\Phi)+(u^r\otimes F^r,\nabla \Phi)+(\nabla u^rF^r,\Phi).
\end{equation}
We rewrite the last term on the right hand side of the latter identity using $\dvr F^r=0$ in $Q_T$ and $u^r=0$ on $\partial\Omega$
\begin{equation}\label{NablaTerm}
	\begin{split}
	\int_0^t(\nabla u^rF^r,\Phi)&=\int_0^t\int_\Omega \partial_k u^r_iF^r_{kj}\Phi_{ij}=\int_0^t\int_\Omega \partial_k (u^r_iF^r_{kj})\Phi_{ij}=-\int_0^t\int_\Omega u^r_iF^r_{kj}\partial_k\Phi_{ij}\\
	&=-\int_0^t(u^r\otimes (F^r)^\top,\nabla\Phi).
	\end{split}
\end{equation}
Multiplying \eqref{RTransport} by $\psi\in C_c(0,T)$, integrating the result over $(0,T)$ and using \eqref{NConvergences}$_{3,5}$, we arrive at
\begin{equation*}
	\int_0^T\psi(t)\bigl((F(t),\Phi(t))-(F_0,\Phi(0))\bigr)=\int_0^T\psi(t)\Biggl(\int_0^t(F,\tder\Phi)+(u\otimes F,\nabla \Phi)-(u\otimes F^\top,\nabla\Phi)\Biggr).
\end{equation*}
Then \eqref{WeakForms}$_2$ follows by the fundamental theorem of the calculus of variations and by the calculations in \eqref{NablaTerm} taking into account that the distributional divergence of $F$ vanishes in $Q_T$. Indeed, assuming that $\{u^\delta\}\subset L^2(0,T;C^\infty_c(\Omega)^d)$ is such that $u^\delta\to u$ in $L^2(0,T;W^{1,2}(\Omega)^d)$ as $\delta\to 0_+$, we have
\begin{equation} \label{analyticalreason}
	\begin{split}
	\int_0^t(u\otimes F^\top,\nabla\Phi)=&\lim_{\delta\to 0_+}\int_0^t\int_\Omega u^\delta_iF_{kj}\partial_k\Phi_{ij}=\lim_{\delta\to 0_+}\int_0^t\int_\Omega F_{kj}\bigl(\partial_k(u^\delta_i\Phi_{ij})-\partial_ku^\delta_i\Phi_{ij}\bigr)\\
	=&-\int_0^t(\nabla uF,\Phi),
	\end{split}
\end{equation}
cf.\ the proof of \cite[Theorem 2.1]{K19} for more details. \colorstart For the last equality we used $\dvr F=0$. \color{black}
In order to obtain \eqref{WeakForms}$_3$ we use $(M^r\times \xi)(t)$ with a fixed $\xi\in C^1(\overline{Q_T})^3$ as a test function in \eqref{ApproxSysClass}$_3$, which is allowed due to the regularity of $M^r$. Integrating over $(0,T)$ we obtain
\begin{equation*}
\begin{split}
	\int_0^T& \bigl(\tder M^r+(u^r\cdot\nabla)M^r\bigr)\cdot M^r\times \xi-M^r\times \bigl(\tder M^r+(u^r\cdot\nabla)M^r-2M^r\times\Hext\bigr)\cdot M^r \times \xi\\&+2\bigl(\nabla M^r, \nabla (M^r\times \xi)\bigr)=0.
	\end{split}
\end{equation*}
We pass to the limit $r\to\infty$ in each term on the left hand side of the latter identity denoting them $I^r_1$, $I^r_2$ and $I^r_3$. Using \eqref{NConvergences}$_{8,11}$, \eqref{HExtReg} and the identities $a\cdot(b\times c)=(a\times b)\cdot c$, $a\times b=-b\times a$ for $a,b,c\in\eR^3$, we get
\begin{equation*}
	\lim_{r\to\infty}I^r_1=-\int_{Q_T} M\times\bigl(\tder M+(u\cdot\nabla)M\bigr) \cdot\xi.
\end{equation*}
Applying the identity $(a\times b)\cdot(c\times d)=(a\cdot c)(b\cdot d)-(a\cdot d)(b\cdot c)$, \eqref{NConvergences}$_{8,9,11}$	and $|M|=1$ a.e.\ in $Q_T$
	\begin{align*}
	\lim_{r\to\infty}I^r_2=&\lim_{r\to\infty}\int_{Q_T} |M^r|^2\bigl(\tder M^r+(u^r\cdot\nabla)M^r-2M^r\times\Hext\bigr)\cdot \xi\\
	&-M^r\cdot\xi\bigl(\tder M^r+(u^r\cdot\nabla)M^r-2M^r\times\Hext\bigr)\cdot M^r\\
	=& \int_{Q_T} \bigl(\tder M+(u\cdot\nabla)M-2M\times\Hext\bigr)\cdot \xi.
\end{align*}
Finally, we integrate by parts and employ \eqref{NConvergences}$_{7,8}$ to get
	\begin{align*}
	\lim_{r\to\infty}I^r_3=&\lim_{r\to\infty}2\int_{Q_T} \sum_{i=1}^d \partial_i M^r\times M^r\cdot\partial_i\xi=\int_{Q_T} 2\sum_{i=1}^d \partial_i M\times M\cdot\partial_i\xi.
\end{align*}
Identity \eqref{WeakForms}$_{3}$ then follows by an integration by parts with respect to time and the density of $C^1(\overline{Q_T})^3$ in $W^{1,2}(Q_T)^3$.

\textbf{Step 4}: Here, we tackle the attainment of the initial data. By the regularity of $\tder u$ provided by Lemma~\ref{Lem:Convergences}, we obtain $u\in C([0,T]; (\VS(\Omega))^*)$, cf.\ \cite[Lemma 7.1]{Ro13}. By the embedding  $L^2_{\dvr}(\Omega)\hookrightarrow (\VS(\Omega))^*$, we thus have that 
\begin{equation}\label{TimeContU}
u\in C_w([0,T];L^2_{\dvr}(\Omega)),
\end{equation}
see \cite[Ch.~III, Lemma~1.4]{Tem77} for details. Next we consider the limit $t\to 0_+$ on both sides of \eqref{WeakForms}$_1$ with $\psi=\psi_1\psi_2$ for an arbitrary but fixed $\psi_1\in C^\infty([0,T])$ such that $\psi_1(0)=1$ and $\psi_2\in C^\infty_c(\Omega)^d$ such that $\dvr\psi_2=0$ in $Q_T$. We then obtain  
\begin{equation*}
	(u(0)-u_0,\psi_2)=0\text{ for all } \psi_2\in C^\infty_c(\Omega)^d,\ \dvr\psi_2=0\text{ in }Q_T,
\end{equation*}
i.e., $u(0)=u_0$ a.e.\ in $\Omega$. First, we infer by similar arguments as above that 
\begin{equation}\label{FMContTime}
\begin{split}
	&F\in C_w([0,T];\CrlH(\Omega)),\\
	&M\in C_w([0,T];W^{1,2}(\Omega)^3)
	\end{split}
\end{equation} 
and 
\begin{equation}\label{FMInVal}
	F(0)=F_0,\ M(0)=M_0\text{ a.e.\ in }\Omega.
\end{equation}
Consequently, we deduce from \eqref{DSEI} that 
\begin{equation*}
\begin{split}
	&\frac{1}{2}\left(\|u(t)\|^2_{L^2(\Omega)}+\|F(t)\|^2_{L^2(\Omega)}+\|\nabla M(t)\|^2_{L^2(\Omega)}\right)+\int_0^t\left(\|\nabla u\|^2_{L^2(\Omega)}+\frac{1}{2}\|\tder M+(u\cdot\nabla)M\|^2_{L^2(\Omega)}\right)\\
	&\leq\int_0^t(\Hext,\tder M)+\frac{1}{2}\left(\|u_0\|^2_{L^2(\Omega)}+\|F_0\|^2_{L^2(\Omega)}+\|\nabla M_0\|^2_{L^2(\Omega)}\right)
	\end{split}
\end{equation*}
holds for all $t\in(0,T)$. Hence we obtain 
\begin{equation*}
\limsup_{t\to 0_+}\left(\|u(t)\|^2_{L^2(\Omega)}+\|F(t)\|^2_{L^2(\Omega)}+\|\nabla M(t)\|^2_{L^2(\Omega)}\right)\leq\|u_0\|^2_{L^2(\Omega)}+\|F_0\|^2_{L^2(\Omega)}+\|\nabla M_0\|^2_{L^2(\Omega)}.
\end{equation*}
On the other hand we have, by \eqref{TimeContU}, \eqref{FMContTime} and \eqref{FMInVal},
\begin{equation*}
\liminf_{t\to 0_+}\left(\|u(t)\|^2_{L^2(\Omega)}+\|F(t)\|^2_{L^2(\Omega)}+\|\nabla M(t)\|^2_{L^2(\Omega)}\right)\geq\|u_0\|^2_{L^2(\Omega)}+\|F_0\|^2_{L^2(\Omega)}+\|\nabla M_0\|^2_{L^2(\Omega)}.
\end{equation*}
The inequalities above, \eqref{TimeContU} and \eqref{FMContTime} imply
\begin{equation*}
	\lim_{t\to 0_+}\left(\|u(t)-u_0\|^2_{L^2(\Omega)}+\|F(t)-F_0\|^2_{L^2(\Omega)}+\|\nabla M(t)-\nabla M_0\|^2_{L^2(\Omega)}\right)=0,
\end{equation*}
which directly implies \eqref{InitCondAtt}.

\section{Appendix}
The first lemma of the appendix deals with equivalent formulations of the Landau-Lifshitz-Gilbert equation under the assumption of a sufficiently regular solution, cf.\ \cite{CF01}. 
\begin{Lemma}\label{Lem:LLGEquivForms}
	Let $\Omega\subset\eR^d$ be a bounded Lipschitz domain, $T>0$, $u\in L^1(Q_T)^d$ , $\Hext\in L^1(Q_T)^3$, $M\in L^1(0,T;W^{2,1}(\Omega)^3)$ and $\tder M\in L^1(Q_T)^3$ be such that $|M|\equiv 1$ a.e.\ in $Q_T$. Then the following forms of the Landau-Lifshitz-Gilbert equation satisfied by $M$ a.e.\ in $Q_T$ are equivalent:
	\begin{align} \label{LG-forms}
	\begin{split}
	\tder M+(u\cdot\nabla)M&=-M\times(\Delta M+\Hext)-M\times\bigl(M\times(\Delta M+\Hext)\bigr),\\%\label{LLGFForm}\\
	\tder M+(u\cdot\nabla)M&=-M\times(\Delta M+\Hext)+\Delta M+\Hext+M\bigl(|\nabla M|^2-M\cdot\Hext\bigr),\\%\label{LLGSForm}\\
	\tder M+(u\cdot\nabla)M&=-2M\times(\Delta M+\Hext)-M\times\bigl(\tder M+(u\cdot\nabla)M\bigr).%\label{LLGTForm}
	\end{split}
	\end{align}
	\begin{proof}
	The equivalence of \eqref{LG-forms}$_1$ and \eqref{LG-forms}$_2$ follows by the aplication of the identities $a\times(b\times c)=(a\cdot c)b-(a\cdot b)c$ fulfilled by any $a,b,c\in \eR^3$ and $M\cdot\Delta M=-|\nabla M|^2$ that is obtained by taking the Laplacian of both sides of $|M|=1$. Using the fact that $M\times M=0$, we get the equivalence of \eqref{LG-forms}$_1$ and
	\begin{equation*}
		\tder M+(u\cdot\nabla)M=-M\times(\Delta M+\Hext)-M\times\left(-M(|\nabla M|^2-M\cdot\Hext)+M\times(\Delta M+\Hext)\right).
	\end{equation*}
	Employing \eqref{LG-forms}$_2$ in the second term on the right hand side of the latter equality, we conclude the equivalence of \eqref{LG-forms}$_1$ and \eqref{LG-forms}$_3$.
	\end{proof}
\end{Lemma}
The ensuing lemma summarizes several assertions concerning the transport equation for the deformation gradient.
\begin{Lemma}\label{Lem:TranspEqRegData}
	Let $\Omega\subset\eR^d$ be a bounded Lipschitz domain, $T>0$, $F_0\in W^{2,\infty}(\Omega)^{d\times d}$, $u\in L^2(0,T;W^{2,\infty}(\Omega)^d)$ with 
	\begin{equation}\label{UConstr}
		\dvr u=0\text{ in }Q_T\text{ and }u=0\text{ on }(0,T)\times\partial\Omega.
	\end{equation}
	Let $F\in L^\infty(0,T;L^2(\Omega)^{d\times d})$ satisfy the transport equation in the sense
	\begin{equation}\label{TrEqWeakly}
		(F(t),\Phi(t))-(F_0,\Phi(0))=\int_0^t(F,\tder\Phi)+(u\otimes F,\nabla\Phi)+(\nabla uF,\Phi)
	\end{equation}
	for all $\Phi\in C^1(\overline{Q_T})^{d\times d}$. Then $F\in C([0,T];W^{2,\infty}(\Omega)^{d\times d})$ with $\tder F\in L^1(0,T;W^{1,\infty}(\Omega)^{d\times d}))$ and $F$ satisfies
		\begin{equation*}
			\tder F+(u\cdot\nabla)F-\nabla uF=0\text{ a.e.\ in }Q_T,\ F(0)=F_0.
		\end{equation*}	
		If the initial datum $F_0$ fulfills additionally $\dvr F_0=0$ in $\Omega$ then $\dvr F=0$ a.e.\ in $Q_T$.
	\begin{proof}
		Let us denote by $\tilde{F}$ a solution of the system of transport equations
		\begin{equation}\label{TFEq}
			\tder \tilde F+(u\cdot\nabla)\tilde F -\nabla u\tilde F=0\text{ in }Q_T,\ \tilde F(0)=F_0
		\end{equation}
		for the initial condition $F_0$ and the velocity $u$ possessing the regularity expressed in the assumptions of the lemma. Then by a standard procedure based on the application of characteristics and the Banach fixed point theorem we get the existence of a unique $ \tilde F\in C([0,T];W^{2,\infty}(\Omega)^{d\times d})$ with $\tder \tilde F\in L^1(0,T;W^{1,\infty}(\Omega)^{d\times d}))$ that satisfies \eqref{TFEq}. Our task is to show that $\tilde F$ coincides with the function $F$ from the assumptions of the lemma. We assume without loss of generality that $F=0$ in $(0,T)\times(\eR^d\setminus \Omega)$ and equation \eqref{TrEqWeakly} being extended from $\Omega$ to $\eR^d$, which is done by adopting ideas from the proof of \cite[Lemma 6.8]{NovStr}. Then we define $F^\varepsilon$ as $F^\varepsilon_{ij}=F_{ij}*\rho_\varepsilon$, $i,j=1,\ldots,d$, where $\rho_\varepsilon$ is defined in \eqref{MolifDef}. The function $F^\varepsilon$ satisfies
		\begin{equation}\label{FEpsEq}
			\tder F^\varepsilon +(u\cdot\nabla)F^\varepsilon-(\nabla uF)^\varepsilon=r^\varepsilon\text{ in }Q_T,\ F^\varepsilon(0)=(F_0)^\varepsilon,
		\end{equation}
		which follows by taking $\Phi^\varepsilon_{ij}=\Phi_{ij}*\rho^\varepsilon$ in the extended version of \eqref{TrEqWeakly} firstly with $\Phi\in C^1_c((0,T)\times\overline\Omega)^{d\times d}$ and secondly with $\Phi\in C^1(\overline{Q_T})^{d\times d}$. We note that $r^\varepsilon\to 0$ as $\varepsilon\to 0$ in $L^2(Q_T)^{d\times d}$, see \cite[Lemma 6.7]{NovStr}. Taking the difference of \eqref{TFEq} and \eqref{FEpsEq}, multiplying the resulting identity with $\tilde F-F^\varepsilon$ and integrating over $Q_t$ with an arbitrary $t\in(0,T]$, applying the solenoidality of $u$ and the fact that $u=0$ on $(0,T)\times \partial\Omega$, we obtain
		\begin{equation}\label{DiffMol}
		\begin{aligned}
			\frac{1}{2}\bigl(\|(\tilde F-F^\varepsilon)(t)\|^2_{L^2(\Omega)}-\|(\tilde F-F^\varepsilon)(0)\|^2_{L^2(\Omega)}\bigr)=&\int_0^t\big(-(u\cdot\nabla)(\tilde F-F^\varepsilon)-(\nabla uF)^\varepsilon\\
			&+\nabla u\tilde F-r^\varepsilon,\tilde F-F^\varepsilon\bigr)\\
			=&\int_0^t\bigl((-\nabla uF)^\varepsilon+\nabla u\tilde{F}-r^\varepsilon,\tilde F-F^\varepsilon\bigr).
		\end{aligned}
		\end{equation}
		Using the facts that 
		\begin{alignat*}{2}
			F^\varepsilon(0)&\to F_0 &&\text{ in }L^2(\Omega)^{d\times d},\\
			F^\varepsilon(t)&\to F(t) &&\text{ in }L^2(\Omega)^{d\times d}\text{ for a.a.\ }t\in(0,T),\\
			(\nabla uF)^\varepsilon(t)&\to (\nabla uF)(t) &&\text{ in }L^2(\Omega)^{d\times d}\text{ for a.a.\ }t\in(0,T)
		\end{alignat*}
		as $\varepsilon\to 0$, we pass to the limit in \eqref{DiffMol} to obtain
		\begin{equation*}
			\frac{1}{2}\|(\tilde F-F)(t)\|^2_{L^2(\Omega)}=\int_0^t\big(\nabla u(\tilde F-F),\tilde F-F\bigr).
		\end{equation*}
		Hence we conclude 
		\begin{equation*}
			\frac{1}{2}\|(\tilde F-F)(t)\|^2_{L^2(\Omega)}\leq \int_0^t\|\nabla u\|_{L^\infty(\Omega)}\|\tilde F- F\|^2_{L^2(\Omega)}.
		\end{equation*}
		Consequently, by the Gronwall lemma we have $F=\tilde F$. To conclude the proof of the lemma, we take the divergence of \eqref{TFEq}, which is allowed due to the amount of regularity possessed by $F$, $\tder F$ and $u$ and obtain
		\begin{equation*}
			\tder \partial_i F_{ij}+\partial_i u_k\partial_k F_{ij}+u_k\partial_k\partial_i F_{ij}=\partial_k\partial_i u_iF_{kj}+\partial_ku_i\partial_iF_{kj}.
		\end{equation*}
		By the solenoidality of $u$ and the switch of indices in the second term on the right hand side, the latter identity becomes
			\begin{equation*}
			\tder \dvr F+u\cdot\nabla \dvr F=0,
		\end{equation*}	
		which is equipped with the initial condition $(\dvr F)(0)=\dvr F_0$ and the lemma is completely proved.
	\end{proof}
\end{Lemma}

\vspace{1em}
\noindent\textbf{Acknowledgement}\\
We are grateful to {\v S}{\'a}rka Ne{\v c}asov{\'a} for her valuable comments on an earlier version of this article and acknowledge financial support by the Deutsche Forschungsgemeinschaft (DFG, German Research Foundation, grant SCHL 1706/4-1, project number 391682204).

\providecommand{\bysame}{\leavevmode\hbox to3em{\hrulefill}\thinspace}
\providecommand{\MR}{\relax\ifhmode\unskip\space\fi MR }
\providecommand{\MRhref}[2]{%
  \href{http://www.ams.org/mathscinet-getitem?mr=#1}{#2}
}
\providecommand{\href}[2]{#2}

\end{document}